\documentclass[12pt]{amsart}
\usepackage{amssymb}
\usepackage[all]{xy}
\usepackage{tikz}
\usepackage{amsmath}
\numberwithin{equation}{section}
\input xy
\xyoption{all}

\newtheorem{lemma}{Lemma}[section]
\newtheorem{corollary}[lemma]{Corollary}
\newtheorem{theorem}[lemma]{Theorem}
\newtheorem{proposition}[lemma]{Proposition}

\theoremstyle{definition}
\newtheorem{remark}[lemma]{Remark}
\newtheorem{definition}[lemma]{Definition}

\DeclareMathOperator{\Mod}{Mod}
\DeclareMathOperator{\modd}{mod}

\DeclareMathOperator{\add}{add}
\DeclareMathOperator{\Hom}{Hom}
\DeclareMathOperator{\pd}{pd}

\DeclareMathOperator{\gl.dim}{gl.dim}
\DeclareMathOperator{\dom.dim}{dom.dim}
\DeclareMathOperator{\Ext}{Ext}

\DeclareMathOperator{\Ker}{Ker}
\DeclareMathOperator{\Coker}{Coker}
\DeclareMathOperator{\Imm}{Im}
\DeclareMathOperator{\E}{E}
\DeclareMathOperator{\eff}{eff}
\DeclareMathOperator{\Eff}{Eff}
\DeclareMathOperator{\gen}{gen}
\DeclareMathOperator{\cogen}{cogen}
\DeclareMathOperator{\adm}{adm}
\DeclareMathOperator{\proj}{proj}

\newtheorem*{theorem a*}{Theorem A}
\newtheorem*{theorem b*}{Theorem B}

\newcounter{diagram}
\numberwithin{diagram}{section}

\newenvironment{diagram}
  {\stepcounter{diagram}\par\smallskip\noindent\begin{minipage}{\linewidth}\centering}
  {\par Diagram~\thediagram\end{minipage}\par\smallskip}

\begin{document}

\title{Higher Auslander correspondence for exact categories}

\author{Ramin Ebrahimi}
\address{Department of Pure Mathematics\\
Faculty of Mathematics and Statistics\\
University of Isfahan\\
P.O. Box: 81746-73441, Isfahan, Iran\\ and School of Mathematics, Institute for Research in Fundamental Sciences (IPM), P.O. Box: 19395-5746, Tehran, Iran}
\email{ramin69@sci.ui.ac.ir / ramin.ebrahimi1369@gmail.com}

\author{Alireza Nasr-Isfahani}
\address{Department of Pure Mathematics\\
Faculty of Mathematics and Statistics\\
University of Isfahan\\
P.O. Box: 81746-73441, Isfahan, Iran\\ and School of Mathematics, Institute for Research in Fundamental Sciences (IPM), P.O. Box: 19395-5746, Tehran, Iran}
\email{nasr$_{-}$a@sci.ui.ac.ir / nasr@ipm.ir}



\begin{abstract}
Inspired by the recent work of Henrard, Kvamme and van Roosmalen \cite{HRK}, we prove a categorified version of higher Auslander correspondence in the context of exact categories. We define $n$-Auslander exact categories and show that there is a bijection between the equivalence classes of $n$-cluster tilting subcategories of exact categories and the equivalence classes of $n$-Auslander exact categories.
\end{abstract}

\maketitle


\section{Introduction}
An artin algebra $\Lambda$ is called
{\em representation-finite} or of {\em finite representation type} provided $\Lambda$ has, up to isomorphism, only finitely many finitely generated indecomposable right $\Lambda$-modules. Artin algebras of finite representation type are of particular importance in representation theory since in this case one has a complete combinatorial description of the module category in terms of the
Auslander--Reiten quiver. In 1971, Auslander defined an artin algebra with nice homological properties which is called {\em Auslander algebra}. In his fundamental theorem, which is called Auslander correspondence, he proved that there is a bijection between Morita equivalence classes of Artin algebras of finite representation type and Auslander algebras \cite{ausla}. We recall that the {\em dominant dimension $\dom.dim \Gamma$} of an artin algebra $\Gamma$ is the largest $m\in \mathbb{N} \cup \{\infty\}$ such that any finitely generated projective right $\Gamma$-module has an injective resolution whose first $m$ terms are projective. An artin algebra $\Gamma$ is called {\em Auslander algebra} if $\gl.dim \Gamma \leq 2\leq\dom.dim \Gamma$.

Beligiannis proved a categorified version of Auslander’s correspondence for abelian categories in \cite{Bel1}. Recall that an abelian category $\mathcal{A}$ with enough projectives and injectives is called an {\em Auslander category} if $\gl.dim \mathcal{A} \leq 2\leq\dom.dim \mathcal{A}$.
Beligiannis in \cite[Theorem 6.6]{Bel1} proved that there is a bijections between equivalence classes of additive categories with split idempotents and abelian Auslander categories.

In representation theory of algebras there are important subcategories of the category of finitely generated right modules which are exact in the
sense of Quillen \cite{Qu}, but not necessarily abelian. For example, the subcategory of the Gorenstein projective modules \cite{AR1, Bel2} and the subcategory of the modules of finite projective dimension \cite{AR, HU}. Recently, Henrard, Kvamme and van Roosmalen in \cite{HRK} defined Auslander exact categories as follows. Let $\mathcal{E}$ be an exact category with enough projectives $\mathcal{P}=\proj(\mathcal{E})$, $^{\bot}\mathcal{P}:=\{E\in\mathcal{E}|\Hom_{\mathcal{E}}(E, P)=0, \forall P\in \mathcal{P}\}$ and $\cogen(\mathcal{P}):=\{E\in\mathcal{E}|
\exists \text{\space inflation \space} E\rightarrowtail P \text{\space with\space} P\in \mathcal{P}\}$. $\mathcal{E}$ is called an {\em Auslander exact category} if it satisfies the following conditions:
\begin{itemize}
\item[(1)]
$(^{\bot}\mathcal{P},\cogen(\mathcal{P}))$ is a torsion pair in $\mathcal{E}$.
\item[(2)]
All morphisms $f:E\rightarrow E'$ with $E'\in ^{\bot}\mathcal{P}$ are admissible with image in $^{\bot}\mathcal{P}$.
\item[(3)]
$\Ext_{\mathcal{E}}^{1}(^{\bot}\mathcal{P},\mathcal{P})=0$.
\item[(4)]
$\gl.dim(\mathcal{E})\leq 2$.
\end{itemize}

Henrard, Kvamme and van Roosmalen proved that there is an equivalence from the 2-category of exact
categories and left exact functors to the 2-category of Auslander exact categories and exact functors
preserving projective objects \cite[Theorem 4.7]{HRK}.

For a positive integer $n$, $n$-Auslander algebras were introduced by Iyama \cite{I2, I1} to construct a higher Auslander correspondence and develop the higher-dimensional analogue of Auslander--Reiten theory. An artin algebra $\Gamma$ is called {\em $n$-Auslander algebra} if $\gl.dim \Gamma \leq n+1\leq\dom.dim \Gamma$. Recall that a finitely generated module $M$ over an artin algebra $\Lambda$ is called an $n$-cluster tilting module if $\{X\in \modd\text{-}\Lambda\mid\Ext^k_\Lambda(M, X)=0, 1\leq k\leq n-1\}=\add M=\{X\in \modd\text{-}\Lambda \mid \Ext^k_\Lambda(X, M)=0, 1\leq k\leq n-1\}$, where $\add M$ denotes the full subcategory of $\modd\text{-}\Lambda$ consisting of the direct factors of finite direct sums of copies of $M$. For more details see \cite{I2}, where the terminology maximal $n$-orthogonal module is used. Iyama in his famous theorem, which is called higher Auslander correspondence, proved that there is a bijection between equivalence classes of $n$-cluster tilting $\Lambda$-modules $M$ for finite-dimensional algebras $\Lambda$ and Morita-equivalence classes of finite-dimensional $n$-Auslander algebras \cite[Theorem, Section 4.2]{I1}.

Beligiannis in \cite{Bel} proved a categorified version of higher Auslander’s correspondence for abelian categories. An abelian category $\mathcal{A}$ is called an {\em $n$-Auslander category} if $\mathcal{A}$ has enough projectives and satisfies $\gl.dim \mathcal{A} \leq n+1\leq\dom.dim \mathcal{A}$. A full subcategory $\mathcal{M}$ of an abelian category $\mathcal{A}$ is called an {\em $n$-cluster tilting subcategory}, if $\mathcal{M}$ is contravariantly finite in $\mathcal{A}$, any right $\mathcal{M}$-approximation is an epimorphism, and $\mathcal{M}$ coincides with the full subcategory $\mathcal{M}^{{\bot}_n}:=\{A\in \mathcal{A}| \Ext^k_\mathcal{A}(\mathcal{M}, A) =0, 1\leq k\leq n-1\}$. Beligiannis proved that there is a bijection between equivalence classes of $n$-cluster tilting subcategories $\mathcal{M}$ in abelian categories $\mathcal{A}$ with enough injectives and equivalence classes of $n$-Auslander categories \cite[Theorem 8.23]{Bel}.

It is natural to ask about a categorified version of higher Auslander’s correspondence for exact categories. Let $\mathcal{E}$ be an exact category and $n$ be a positive integer. According to \cite[Definition 4.13]{J},
a subcategory $\mathcal{M}$ of $\mathcal{E}$ is called {\em $n$-cluster tilting subcategory} if the following conditions are satisfied.
\begin{itemize}
\item[$(i)$]
Every object $E\in \mathcal{E}$ has a left $\mathcal{M}$-approximation by an inflation $E\rightarrowtail M.$
\item[$(ii)$]
Every object $E\in \mathcal{E}$ has a right $\mathcal{M}$-approximation by a deflation $M'\twoheadrightarrow E.$
\item[$(iii)$]
We have
\begin{align}
\mathcal{M}&=\{ E\in \mathcal{E} \mid \forall i\in \{1, \ldots, n-1 \}, \Ext_{\mathcal{E}}^i(\mathcal{M},E)=0 \} \notag\\
&= \{ E\in \mathcal{E} \mid \forall i\in \{1, \ldots, n-1 \}, \Ext_{\mathcal{E}}^i(E,\mathcal{M})=0 \}.\notag
\end{align}
\end{itemize}

Let $\mathcal{E}$ be an exact category with enough projectives $\mathcal{P}=\proj(\mathcal{E})$ and $n$ be a positive integer. We say that $\mathcal{E}$ is an {\em $n$-Auslander exact category} if it satisfies the following conditions:
\begin{itemize}
\item[(1)]
$(^{\bot}\mathcal{P},\cogen(\mathcal{P}))$ is a torsion pair in $\mathcal{E}$.
\item[(2)]
All morphisms $f:E\rightarrow E'$ with $E'\in ^{\bot}\mathcal{P}$ are admissible with image in $^{\bot}\mathcal{P}$.
\item[(3)]
$\Ext_{\mathcal{E}}^{1,\cdots,n}(^{\bot}\mathcal{P},\mathcal{P})=0$.
\item[(4)]
$\gl.dim(\mathcal{E})\leq n+1$.
\item[(5)]
$\mathcal{P}$ is an admissibly covariantly finite subcategory of $\mathcal{E}$.
\end{itemize}

Note that, by \cite[Lemma 4.16]{HRK}, $1$-Auslander exact categories are precisely Auslander exact categories in sense of Henrard, Kvamme and van Roosmalen \cite[Definition 4.1]{HRK}.
In this paper we prove a categorified version of higher Auslander’s correspondence for exact categories. More precisely we prove the following theorem.

\begin{theorem a*}$($Theorem \ref{th2}$)$
There is a bijection between the following:
\begin{itemize}
\item[(1)]
Equivalence classes of $n$-cluster tilting subcategories of exact categories.
\item[(2)]
Equivalence classes of $n$-Auslander exact categories.
\end{itemize}
\end{theorem a*}

We also prove the more familiar version of higher Auslander correspondence for exact categories.

\begin{theorem b*}$($Theorem \ref{th3}$)$
There is a bijection between the following:
\begin{itemize}
\item[(1)]
Equivalence classes of $n$-cluster tilting subcategories of exact categories with enough injectives.
\item[(2)]
Equivalence classes of exact categories $\mathcal{E}$ with enough projectives $\mathcal{P}=\proj(\mathcal{E})$ satisfying the following conditions:

$(a)$
$\gl.dim(\mathcal{E})\leq n+1\leq \dom.dim(\mathcal{E})$.

$(b)$
Any morphism $X\rightarrow E$ with $E\in ^{\bot}\mathcal{P}$ is admissible.

$(c)$
$\mathcal{P}$ is admissibly covariantly finite.
\end{itemize}
\end{theorem b*}

The paper is organized as follows. In Section 2 we recall the definition and some basic properties of exact categories. We also recall the definition and some properties of n-cluster tilting subcategories of exact categories and prove a criterion for n-cluster tilting subcategories of exact categories that will be use in the proof of our main theorem, see Proposition \ref{pro15}. We end by recalling the definition and properties of localisation in exact categories. In Section 3 we recall the definition of admissibly presented functors and show some properties of the category of admissibly presented functors $\modd_{\adm}\mathcal{E}$. When $\mathcal{E}$ is an exact category and $\mathcal{M}$ is an $n$-cluster tilting subcategory of it we show that there is an adjoint pair $(\mathbb{R},\mathbb{E}):\modd_{\adm}(\mathcal{E})\rightarrow \modd_{\adm}(\mathcal{M})$ that is a restriction of the adjoint pair $(\mathbb{R},\mathbb{E}):\Mod(\mathcal{E})\rightarrow \Mod(\mathcal{M})$. We also define a subcategory $\mathcal{F}(\mathcal{M})$ of $\modd_{\adm}\mathcal{M}$ and show that $(\eff(\mathcal{M}),\mathcal{F}(\mathcal{M}))$ is a torsion pair in $\modd_{\adm}\mathcal{M}$. Finally we study some properties of this torsion pair. In Section 4 we first define $n$-Auslander exact categories and then we show that if $\mathcal{E}$ is an $n$-Auslander exact category with $\mathcal{P}=\proj(\mathcal{E})$, then the composition $H:\mathcal{P}\hookrightarrow \mathcal{E}\rightarrow \frac{\mathcal{E}}{^{\bot}\mathcal{P}}$ is fully faithful and it's essential image is an $n$-cluster tilting subcategory of $\frac{\mathcal{E}}{^{\bot}\mathcal{P}}$. By using this result we prove Theorem \ref{th2}. Finally, in Section 5 we prove Theorem \ref{th3}.

\subsection{Notation}
Throughout this paper, all categories are assumed to be essentially small and additive. All subcategories are
assumed to be full and closed under isomorphisms.


\section{preliminaries}
In this section we recall the definition of exact category, which is introduced by Quillen in \cite{Qu}, and give some result that will be needed later. We follow closely \cite{Bu}.

Let $\mathcal{E}$ be an additive category. A {\em kernel-cokernel pair} in $\mathcal{E}$ is a pair
of composable morphisms
\begin{equation}
X\overset{i}{\rightarrow}Y\overset{p}{\rightarrow}Z\notag
\end{equation}
such that $i$ is a kernel of $p$ and $p$ is a cokernel of $i$. A {\em conflation category} is an additive category with a class of kernel-cokernel pairs, called {\em conflations}, closed under isomorphisms. A map that occurs as the kernel (resp., cokernel) in a conflation is called an {\em inflation} (resp., {\em deflation}). As usual we depict inflations by $\rightarrowtail$
and deflations by $\twoheadrightarrow$ \cite[Definition ]{Bu}. A map $f : X \rightarrow Y$ is called an {\em admissible} if admits a factorization $X\twoheadrightarrow X'\rightarrowtail Y$ \cite[Definition 2.1]{HRK}.

\begin{definition}$($\cite[Definition ]{Bu}$)$
Let $\mathcal{E}$ be an additive category. An {\em exact structure}
on $\mathcal{E}$ is a class $\mathcal{X}$ of kernel-cokernel pairs in $\mathcal{E}$, closed under isomorphisms which satisfies the following axioms:
\begin{itemize}
\item[$(\E0)$]
For all objects $E\in \mathcal{E}$, the identity morphism $1_E$ is an inflation.
\item[$(\E0^{op})$]
For all objects $E\in \mathcal{E}$, the identity morphism $1_E$ is a deflation.
\item[$(\E1)$]
The class of inflations is closed under composition.
\item[$(\E1^{op})$]
The class of deflations is closed under composition.
\item[$(\E2)$]
The pushout of an inflation along an arbitrary morphism exists and yields
an inflation. The situation is illustrated in the following commutative diagram:
\begin{center}
\begin{tikzpicture}
\node (X1) at (-4,1) {$X^0$};
\node (X2) at (-2,1) {$X^1$};
\node (X6) at (-4,-1) {$Y^0$};
\node (X7) at (-2,-1) {$Y^1$};
\draw [>->,thick] (X1) -- (X2) node [midway,above] {};
\draw [>->,thick,dashed] (X6) -- (X7) node [midway,above] {};
\draw [->,thick] (X1) -- (X6) node [midway,left] {};
\draw [->,thick,dashed] (X2) -- (X7) node [midway,left] {};
\end{tikzpicture}
\end{center}

\item[$(\E2^{op})$]
The pullback of a deflation along an arbitrary morphism exists and yields
a deflation. The situation is illustrated in the following commutative diagram:
\begin{center}
\begin{tikzpicture}
\node (X4) at (2,1) {$X^n$};
\node (X5) at (4,1) {$X^{n+1}$};
\node (X9) at (2,-1) {$Y^n$};
\node (X10) at (4,-1) {$Y^{n+1}$};
\draw [->>,thick,dashed] (X4) -- (X5) node [midway,above] {};
\draw [->>,thick] (X9) -- (X10) node [midway,above] {};
\draw [->,thick] (X5) -- (X10) node [midway,right] {};
\draw [->,thick,dashed] (X4) -- (X9) node [midway,left] {};
\end{tikzpicture}
\end{center}
\end{itemize}
\end{definition}

An {\em exact category} is a pair $(\mathcal{E},\mathcal{X})$ where $\mathcal{E}$ is an additive category and $\mathcal{X}$ is an exact structure on $\mathcal{E}$. If the class $\mathcal{X}$ is clear from the context, we identify $\mathcal{E}$ with the
pair $(\mathcal{E},\mathcal{X})$. The members of $\mathcal{X}$ are called {\em conflations} or {\em admissible short exact sequences}.

\begin{remark}

\begin{itemize}
\item[$(1)$]
An exact category in the sense of Quillen (see \cite{Qu}) is a conflation category $\mathcal{E}$
satisfying axioms $\E0$ through $\E3$ and $\E0^{op}$ through $\E3^{op}$, where
\begin{itemize}
\item[$(\E3)$]
If $f:A\rightarrow B$ and $g: B\rightarrow C$ are morphisms in $\mathcal{E}$ such that $f$ has a cokernel and $gf$ is an
inflation, then $f$ is an inflation.
\item[$(\E3^{op})$]
If $f:A\rightarrow B$ and $g: B\rightarrow C$ are morphisms in $\mathcal{E}$ such that $g$ has a kernel and $gf$ is a deflation,
then $g$ is a deflation.
\end{itemize}
\item[$(2)$] Keller, in \cite[appendix A]{Kl}, shows that
axioms $(\E0^{op})$, $(\E1^{op})$, $(\E2^{op})$, and $(\E2)$ suffice to define an exact category.
\item[$(3)$] Axioms $(\E3)$ and $(\E3^{op})$ are sometimes referred to as obscure axioms (see \cite{Bu}).
\end{itemize}
\end{remark}

The following proposition will be used frequently in the rest of the paper.

\begin{proposition}\label{pro3}
Let $\mathcal{E}$ be an exact category.
\begin{itemize}
\item[(1)]
Let $g:Y\rightarrow Z$ be a morphism in $\mathcal{E}$ that admit a kernel. If there is a morphism $g':Y'\rightarrow Y$ such that $gg'$ is a deflation, then $g$ is a deflation.
\item[(2)]
Let $g:Y\rightarrow Z$ be a morphism in $\mathcal{E}$. If there is a morphism $g':Y'\rightarrow Y$ such that $gg'$ is a deflation, then $(0,g):Y'\oplus Y\rightarrow Z$ is a deflation.
\end{itemize}
\begin{proof}
$(1)$ is the obscure axiom \cite[Proposition 2.16]{Bu}. For $(2)$ see the \cite[Proposition 3.4]{HR2} and it's proof.
\end{proof}
\end{proposition}

A cochain $\cdots\overset{d^{n-2}}{\longrightarrow} X^{n-1}\overset{d^{n-1}}{\longrightarrow} X^{n}\overset{d^{n}}{\longrightarrow} X^{n+1}\overset{d^{n+1}}{\longrightarrow}\cdots$ in a conflation category $\mathcal{C}$ is called {\em acyclic} or
{\em exact} if each $d^i$ is admissible and $\Ker(d^{i+1})=\Imm(d^i)$ \cite[Definition 2.1]{HRK}.

\begin{definition}$($\cite[Definition 2.13]{HRK}$)$
Let $(\mathcal{E},\mathcal{X})$ and $(\mathcal{E}',\mathcal{X}')$ be exact categories, and $F:\mathcal{E}\rightarrow \mathcal{E}'$ a covariant additive functor.
\begin{itemize}
\item[(1)]
$F$ is called {\em left exact} if for every conflation $X\overset{f}{\rightarrowtail} Y\overset{g}{\twoheadrightarrow} Z$ in $(\mathcal{E},\mathcal{X})$, $F(X)\overset{F(f)}{\rightarrowtail} F(Y)\overset{F(g)}{\rightarrow} F(Z)\twoheadrightarrow \Coker(F(g))$ be an acyclic complex in $(\mathcal{E}',\mathcal{X}')$.
\item[(2)]
$F$ is called {\em right exact} if for every conflation $X\overset{f}{\rightarrowtail} Y\overset{g}{\twoheadrightarrow} Z$ in $(\mathcal{E},\mathcal{X})$, $\Ker(F(g))\rightarrowtail F(X)\overset{F(f)}{\rightarrow} F(Y)\overset{F(g)}{\twoheadrightarrow} F(Z)$ be an acyclic complex in $(\mathcal{E}',\mathcal{X}')$.
\item[(3)]
$F$ is called {\em exact} if for every conflation $X\overset{f}{\rightarrowtail} Y\overset{g}{\twoheadrightarrow} Z$ in $(\mathcal{E},\mathcal{X})$, $F(X)\overset{F(f)}{\rightarrowtail} F(Y)\overset{F(g)}{\twoheadrightarrow} F(Z)$ belongs to $\mathcal{X}'$.
\end{itemize}
\end{definition}

The contravariant left exact (resp., right exact, exact) functors are defined similarly.

Let $\mathcal{M}$ be an additive category and $f:A\rightarrow B$ a morphism in $\mathcal{M}$. A {\em weak cokernel} of $f$ is a morphism $g:B\rightarrow C$ such that for all $C^{\prime} \in \mathcal{M}$  the sequence of abelian groups
\begin{equation}
\Hom_\mathcal{M}(C,C')\overset{\Hom_\mathcal{M}(g,C')}{\longrightarrow} \Hom_\mathcal{M}(B,C')\overset{\Hom_\mathcal{M}(f,C')}{\longrightarrow} \Hom_\mathcal{M}(A,C') \notag
\end{equation}
is exact. The concept of {\em weak kernel} is defined dually.

Let $d^0:X^0 \rightarrow X^1$ be a morphism in $\mathcal{M}$. An {\em $n$-cokernel} of $d^0$ is a sequence
\begin{equation}
(d^1, \ldots, d^n): X^1 \overset{d^1}{\rightarrow} X^2 \overset{d^2}{\rightarrow}\cdots \overset{d^{n-1}}{\rightarrow} X^n \overset{d^n}{\rightarrow} X^{n+1} \notag
\end{equation}
of objects and morphisms in $\mathcal{M}$ such that for each $Y\in \mathcal{M}$
the induced sequence of abelian groups
\begin{align}
0 \rightarrow \Hom_\mathcal{M}(X^{n+1},Y) \rightarrow \Hom_\mathcal{M}(X^n,Y) \rightarrow\cdots\rightarrow \Hom_\mathcal{M}(X^1,Y) \rightarrow \Hom_\mathcal{M}(X^0,Y) \notag
\end{align}
is exact. Equivalently, the sequence $(d^1, \ldots, d^n)$ is an $n$-cokernel of $d^0$ if for all $1\leq k\leq n-1$
the morphism $d^k$ is a weak cokernel of $d^{k-1}$ and $d^n$ is moreover a cokernel of $d^{n-1}$ \cite[Definition 2.2]{J}. The concept of {\em $n$-kernel} of a morphism is defined dually.

\begin{definition}(\cite[Definition 2.4]{L})\label{d1}
Let $\mathcal{M}$ be an additive category. A {\em right $n$-exact sequence} in $\mathcal{M}$ is a complex
\begin{equation}
X^0 \overset{d^0}{\rightarrow} X^1 \overset{d^1}{\rightarrow} \cdots \overset{d^{n-1}}{\rightarrow} X^n \overset{d^n}{\rightarrow} X^{n+1} \notag
\end{equation}
such that $(d^1, \ldots, d^{n})$ is an $n$-cokernel of $d^0$. The concept of {\em left $n$-exact sequence} is defined dually. An {\em $n$-exact sequence} is a sequence which is both a right $n$-exact sequence and a left $n$-exact sequence.
\end{definition}

Let $\mathcal{A}$ be an additive category and $\mathcal{B}$ be a full subcategory of $\mathcal{A}$. $\mathcal{B}$ is called {\em covariantly finite in $\mathcal{A}$} if for every $A\in \mathcal{A}$ there exists an object $B\in\mathcal{B}$ and a morphism
$f : A\rightarrow B$ such that, for all $B'\in\mathcal{B}$, the sequence of abelian groups $\Hom_\mathcal{A}(B, B')\rightarrow \Hom_\mathcal{A}(A, B')\rightarrow 0$ is exact. Such a morphism $f$ is called a {\em left $\mathcal{B}$-approximation of $A$}. The notions of {\em contravariantly
finite subcategory of $\mathcal{A}$} and {\em right $\mathcal{B}$-approximation} are defined dually. A {\em functorially
finite subcategory of $\mathcal{A}$} is a subcategory which is both covariantly and contravariantly finite
in $\mathcal{A}$ \cite[page 113]{AR}.

\begin{definition}$($\cite[Definition 4.13]{J}$)$
Let $(\mathcal{E},\mathcal{X})$ be an exact category and $\mathcal{M}$ a subcategory of $\mathcal{E}$. $\mathcal{M}$ is called an $n$-cluster tilting subcategory of $(\mathcal{E},\mathcal{X})$ if the following conditions are satisfied.
\begin{itemize}
\item[$(i)$]
Every object $E\in \mathcal{E}$ has a left $\mathcal{M}$-approximation by an inflation $E\rightarrowtail M$.
\item[$(ii)$]
Every object $E\in \mathcal{E}$ has a right $\mathcal{M}$-approximation by a deflation $M'\twoheadrightarrow E$.
\item[$(iii)$]
We have
\begin{align}
\mathcal{M}&=\{ E\in \mathcal{E} \mid \forall i\in \{1, \ldots, n-1 \}, \Ext_{\mathcal{E}}^i(\mathcal{M},E)=0 \} \notag\\
&= \{ E\in \mathcal{E} \mid \forall i\in \{1, \ldots, n-1 \}, \Ext_{\mathcal{E}}^i(E,\mathcal{M})=0 \}.\notag
\end{align}
\end{itemize}
Note that $\mathcal{E}$ itself is the unique $1$-cluster tilting subcategory of $\mathcal{E}$.
\end{definition}

A subcategory $\mathcal{M}$ of an exact category $\mathcal{E}$ is called {\em generating} if for every object
$X\in \mathcal{E}$ there exists an object $Y\in\mathcal{M}$ and a deflation $Y\twoheadrightarrow X$. The concept of
{\em cogenerating} subcategory is defined dually \cite[Definition 3.13]{HRK}.

\begin{remark}\label{rm1}
Let $\mathcal{E}$ be an exact category and $\mathcal{M}$ be a full subcategory of $\mathcal{E}$. By the proof of \cite[Proposition 4.3]{Kv1}, if $f:M\rightarrow E$ is a right $\mathcal{M}$-approximation and $g:M'\rightarrow E$ is a deflation with $M'\in \mathcal{M}$ then $(f,g):M\oplus M'\rightarrow E$ is a deflation and right $\mathcal{M}$-approximation. This shows that the following are equivalent.
\begin{itemize}
\item[(1)]
Every object $E\in \mathcal{E}$ has a right $\mathcal{M}$-approximation by a deflation $M'\twoheadrightarrow E$.
\item[(2)]
$\mathcal{M}$ is a generating and contravariantly finite subcategory of $\mathcal{E}$.
\end{itemize}
Dually, we can see that the following are equivalent.
\begin{itemize}
\item[(1)]
Every object $E\in \mathcal{E}$ has a left $\mathcal{M}$-approximation by an inflation $E\rightarrowtail M'$.
\item[(2)]
$\mathcal{M}$ is a cogenerating and covariantly finite subcategory of $\mathcal{E}$.
\end{itemize}
\end{remark}

A full subcategory $\mathcal{M}$ of an exact or abelian category $\mathcal{E}$ is called {\em $n$-rigid}, if for every
two objects $M,N\in \mathcal{M}$ and for every $k\in \{1,\cdots, n - 1\}$, we have $\Ext_{\mathcal{E}}^k(M, N)=0$ \cite[Page 443]{Bel}. Any $n$-cluster tilting subcategory $\mathcal{M}$ of an exact category $\mathcal{E}$ is $n$-rigid.

\begin{proposition}\label{pro1}
Let $\mathcal{E}$ be an exact category and $\mathcal{M}$ be an $n$-rigid subcategory of $\mathcal{E}$ that is closed under direct summand.
\begin{itemize}
\item[(1)]
Assume that $X\in \mathcal{E}$ and there exists an exact sequence
\begin{equation}
0\rightarrow M_n\rightarrow M_{n-1}\rightarrow \cdots\rightarrow M_1\rightarrow X\rightarrow 0 \notag
\end{equation}
with $M_i\in \mathcal{M}$ for each $i\in \{1,\cdots,n\}$. Then $M_1\rightarrow X$ is a right $\mathcal{M}$-approximation and if $\Ext_{\mathcal{E}}^k(X,\mathcal{M})=0$ for all $k\in \{1,\cdots,n-1\}$, then $X\in \mathcal{M}$.
\item[(2)]
Dually, Assume that $X\in \mathcal{E}$ and there exists an exact sequence
\begin{equation}
0\rightarrow X\rightarrow M_{1}\rightarrow \cdots\rightarrow M_{n-1}\rightarrow M_{n}\rightarrow 0 \notag
\end{equation}
with $M_i\in \mathcal{M}$ for each $i\in \{1,\cdots,n\}$. Then $X\rightarrow M_{1}$ is a left $\mathcal{M}$-approximation and if $\Ext_{\mathcal{E}}^k(\mathcal{M}, X)=0$ for all $k\in \{1,\cdots,n-1\}$, then $X\in \mathcal{M}$.
\end{itemize}
\begin{proof}
See the proof of \cite[Proposition 4.4]{Kv1}.
\end{proof}
\end{proposition}

The following criterion for $n$-cluster tilting subcategories of exact categories is useful for the proof of our main theorem.

\begin{proposition}\label{pro15}
Let $\mathcal{E}$ be an exact category and $\mathcal{M}$ be a subcategory of $\mathcal{E}$ that is a generating-cogenerating, functorially finite and $n$-rigid subcategory. Then the followings are equivalent.
\begin{itemize}
\item[(1)]
$\mathcal{M}$ is an $n$-cluster tilting subcategory of $\mathcal{E}$.
\item[(2)]
$\mathcal{M}$ is closed under direct summand and for every $X\in\mathcal{E}$ there is an exact sequence
\begin{equation}
0\rightarrow M_n\rightarrow M_{n-1}\rightarrow \cdots\rightarrow M_1\rightarrow X\rightarrow 0 \notag
\end{equation}
in $\mathcal{E}$ with $M_i\in \mathcal{M}$ for each $i\in\{1,\cdots,n\}$.
\end{itemize}
\begin{proof}
The implication $(1)\Rightarrow (2)$ follows from the dual of \cite[Proposition 4.15]{J}. We show the implication $(2)\Rightarrow (1)$. Let $X\in \mathcal{E}$ such that $\Ext_{\mathcal{E}}^{1,\ldots,n-1}(X,\mathcal{M})=0$. By Proposition \ref{pro1} we have that $X\in \mathcal{M}$.

Now let $X\in \mathcal{E}$ such that $\Ext_{\mathcal{E}}^{1,\ldots,n-1}(\mathcal{M},X)=0$. Choose a conflation $X\rightarrowtail M^1\twoheadrightarrow K^1$ where $X\rightarrowtail M^1$ is a left $\mathcal{M}$-approximation. By applying $\Hom_{\mathcal{E}}(-,N)$ with $N\in \mathcal{M}$ to this conflation we see that $\Ext_{\mathcal{E}}^1(K^1,\mathcal{M})=0$. Again choose a conflation $K^1\rightarrowtail M^2\twoheadrightarrow K^2$ where $K^1\rightarrowtail M^2$ is a left $\mathcal{M}$-approximation. Similarly $\Ext_{\mathcal{E}}^1(K^2,\mathcal{M})=0$ and $\Ext_{\mathcal{E}}^2(K^2,\mathcal{M})\cong \Ext_{\mathcal{E}}^1(K^1,\mathcal{M})=0$. Repeating this argument we obtain exact sequence
\begin{equation}
0\rightarrow X\rightarrow M^1\rightarrow M^2\rightarrow \cdots\rightarrow M^{n-1}\rightarrow K^{n-1}\rightarrow 0 \notag
\end{equation}
with $\Ext_{\mathcal{E}}^{1,\ldots,n-1}(K^{n-1},\mathcal{M})=0$. By the first part of the proof we have $K^{n-1}\in \mathcal{M}$. Then the result follows from Proposition \ref{pro1}.
\end{proof}
\end{proposition}

In the following proposition we recall some of the basic properties of $n$-cluster tilting subcategories.

\begin{proposition}
Let $(\mathcal{E},\mathcal{X})$ be an exact category and $\mathcal{M}$ be an $n$-cluster
tilting subcategory of $(\mathcal{E},\mathcal{X})$. Then
\begin{itemize}
\item[(1)]
Every deflation $X^n \overset{d^n}{\twoheadrightarrow} X^{n+1}$ fits into an $n$-exact sequence
\begin{equation}\label{nex}
X^0 \overset{d^0}{\rightarrowtail} X^1 \overset{d^1}{\rightarrow} \cdots \overset{d^{n-1}}{\rightarrow} X^n \overset{d^n}{\twoheadrightarrow} X^{n+1}
\end{equation}
in $\mathcal{M}$, that is an $\mathcal{X}$-acyclic complex.
\item[(2)]
For every $n$-exact sequence of the form \eqref{nex} and each morphism $g:Y^{n+1}\rightarrow X^{n+1}$ there exists commutative diagram
\begin{center}
\begin{tikzpicture}
\node (X1) at (-4,1) {$Y^0$};
\node (X2) at (-2,1) {$Y^1$};
\node (X3) at (0,1) {$\cdots$};
\node (X4) at (2,1) {$Y^n$};
\node (X5) at (4,1) {$Y^{n+1}$};
\node (X6) at (-4,-1) {$X^0$};
\node (X7) at (-2,-1) {$X^1$};
\node (X8) at (0,-1) {$\cdots$};
\node (X9) at (2,-1) {$X^n$};
\node (X10) at (4,-1) {$X^{n+1}$};
\draw [->,thick,dashed] (X1) -- (X2) node [midway,above] {$d_Y^0$};
\draw [->,thick,dashed] (X2) -- (X3) node [midway,above] {$d_Y^1$};
\draw [->,thick,dashed] (X3) -- (X4) node [midway,above] {$d_Y^{n-1}$};
\draw [->>,thick,dashed] (X4) -- (X5) node [midway,above] {$d_Y^n$};
\draw [->,thick] (X6) -- (X7) node [midway,above] {$d_X^0$};
\draw [->,thick] (X7) -- (X8) node [midway,above] {$d_X^1$};
\draw [->,thick] (X8) -- (X9) node [midway,above] {$d_X^{n-1}$};
\draw [->>,thick] (X9) -- (X10) node [midway,above] {$d_X^{n}$};
\draw [->,thick,dashed] (X1) -- (X6) node [midway,left] {};
\draw [->,thick] (X5) -- (X10) node [midway,right] {$g$};
\draw [->,thick,dashed] (X2) -- (X7) node [midway,left] {};
\draw [->,thick,dashed] (X4) -- (X9) node [midway,left] {};
\end{tikzpicture}
\end{center}
where rows are $n$-exact sequences.
\end{itemize}
\begin{proof}
Follows from \cite[Theorem 4.14]{J} and \cite[Proposition 4.15]{J}.
\end{proof}
\end{proposition}

\subsection{Localisation of exact categories}
In this subsection we recall the localisation theory of exact categories.

\begin{definition}$($\cite[Definition 2.4]{HRK}$)$
Let $\mathcal{E}$ be an exact category. A full nonempty subcategory $\mathcal{A}\subseteq\mathcal{E}$ is called an {\em admissibly deflation-percolating subcategory} if the following axioms are satisfied:
\begin{itemize}
\item[\textbf{A1}]
The category $\mathcal{A}$ is a {\em Serre subcategory} of $\mathcal{E}$, i.e. given a conflation $X\rightarrowtail Y\twoheadrightarrow Z$ in $\mathcal{E}$, we have
that $Y\in\mathcal{A} \Leftrightarrow X, Z \in\mathcal{A}.$
\item[\textbf{A2}]
All morphisms $f : X \rightarrow A$ with $A\in \mathcal{A}$ are admissible with image in $\mathcal{A}$, i.e. factors as
\begin{equation}
X\twoheadrightarrow\Imm(f)\rightarrowtail A \notag
\end{equation}
with $\Imm(f)\in \mathcal{A}$.
\end{itemize}
Dually, $\mathcal{A}$ is called an {\em admissibly inflation-percolating subcategory} if $\mathcal{A}^{op}$ is an admissibly deflation-percolating subcategory in $\mathcal{E}^{op}$. Finally, $\mathcal{A}$ is called a {\em two-sided admissibly percolating subcategory} if it is both
an admissibly deflation and admissibly inflation-percolating subcategory.
\end{definition}

\begin{definition}$($\cite[Definition 2.6]{HRK}$)$
Let $\mathcal{E}$ be an exact category and $\mathcal{A}\subseteq\mathcal{E}$ be an admissibly deflation-percolating subcategory. A morphism $f:X\rightarrow Y$ in $\mathcal{E}$ is called a {\em weak $\mathcal{A}$-isomorphism} or simply a {\em weak isomorphism} if $\mathcal{A}$ is understood, if $f$ is admissible and $\Ker(f),\Coker(f)\in\mathcal{A}$. The set of all weak $\mathcal{A}$-isomorphisms is denoted by $\Sigma_{\mathcal{A}}$
\end{definition}

In the following proposition we collect some of the properties of localisation of an exact category with respect to a two-sided admissibly percolating subcategory.

\begin{proposition}$($\cite[Proposition 2.7, Theorem 2.8 and Theorem 2.10]{HRK}$)$\label{pro11}
Let $\mathcal{E}$ be an exact category and $\mathcal{A}\subseteq\mathcal{E}$ be a two-sided admissibly percolating subcategory.
\begin{itemize}
\item[(1)]
The set $\Sigma_{\mathcal{A}}$ is a left and right multiplicative system.
\item[(2)]
Let $Q:\mathcal{E}\rightarrow \Sigma_{\mathcal{A}}^{-1}\mathcal{E}$ be the localisation functor. Then
\begin{itemize}
\item[(i)]
$Q(M)=0$ for an object $M\in \mathcal{E}$ if and only if $M\in \mathcal{A}$.
\item[(ii)]
$Q(f)$ is an isomorphism if and only if $f\in \Sigma_{\mathcal{A}}$.
\item[(iii)]
$f$ is an admissible morphism in $\mathcal{E}$ if and only if $Q(f)$ is admissible morphism in $\Sigma_{\mathcal{A}}^{-1}\mathcal{E}$.
\end{itemize}
\item[(3)]
$Q$ sends conflations of $\mathcal{E}$ to kernel-cokernel pairs in $\Sigma_{\mathcal{A}}^{-1}\mathcal{E}$. Furthermore, closing this class of kernel-cokernel pairs under isomorphism makes $\Sigma_{\mathcal{A}}^{-1}\mathcal{E}$ into an exact category.
\item[(4)]
The functor $Q$ is the universal exact functor from $\mathcal{E}$ to an exact category that annihilate $\mathcal{A}$.
\end{itemize}
\end{proposition}
\begin{definition}$($\cite[Definition 2.11]{HRK}$)$
Let $\mathcal{E}$ be an exact category. A pair $(\mathcal{T},\mathcal{F})$ of subcategories of $\mathcal{E}$ closed under isomorphism is called a {\em torsion pair}, if
\begin{itemize}
\item[(1)]
For every $T\in\mathcal{T}$ and $F\in \mathcal{F}$ we have $\Hom_\mathcal{E}(T,F)=0$;
\item[(2)]
For any object $E\in \mathcal{E}$, there exists a conflation $T\rightarrowtail E\twoheadrightarrow F$ with $T\in\mathcal{T}$ and $F\in \mathcal{F}$.
\end{itemize}
\end{definition}

For properties of torsion pair we refer the reader to \cite[Section 2.4]{HRK}. We need the following proposition in the rest of the paper.

\begin{proposition}$($\cite[Proposition 2.16]{HRK}$)$
Let $\mathcal{E}$ be an exact category with a torsion pair $(\mathcal{T},\mathcal{F})$ such that $\mathcal{T}$ satisfies axiom $\textbf{A2}$. Then the subcategory $\mathcal{T}\subseteq \mathcal{E}$ is two-sided admissibly percolating.
\end{proposition}

\section{The category of admissibly presented functors}
In this section, for an $n$-cluster tilting subcategory $\mathcal{M}$ of an exact category $\mathcal{E}$ we study some properties of admissibly presented functors from $\mathcal{M}$ to the category of all abelian groups.

We denote by $\Mod (\mathcal{M})$, the category of all additive contravariant functors $F:\mathcal{M}\rightarrow Ab$ from $\mathcal{M}$ to the category of abelian groups. Let
\begin{align*}
\Upsilon_{\mathcal{M}}:&\mathcal{E}\longrightarrow \Mod (\mathcal{M})  \notag \\
&X\longmapsto \Hom_{\mathcal{E}}(-,X)|_{\mathcal{M}} \notag
\end{align*}
In particular, $\Upsilon_{\mathcal{E}}$ is the Yoneda embedding $\Upsilon_{\mathcal{E}}:\mathcal{E}\hookrightarrow \Mod (\mathcal{E})$.

\begin{definition}
Let $\mathcal{M}$ be a full subcategory of an exact category $\mathcal{E}$.
\begin{itemize}
\item[(1)]
Let $\modd (\mathcal{M})$ the full subcategory of $\Mod (\mathcal{M})$ consisting of all {\em finitely presented functors}, i.e. those functors $F$ that admit a projective presentation
\begin{equation}
\Upsilon_{\mathcal{M}}(M)\overset{\Upsilon_{\mathcal{M}}(f)}{\longrightarrow} \Upsilon_{\mathcal{M}}(N)\rightarrow F\rightarrow 0,\notag
\end{equation}
for some morphism $f : M\rightarrow N$ in $\mathcal{M}$.
\item[(2)]
Let $\modd_{adm} (\mathcal{M})$ be the full subcategory of $\Mod (\mathcal{M})$ consisting of those functors $F$ that admit a projective presentation
\begin{equation}
\Upsilon_{\mathcal{M}}(M)\overset{\Upsilon_{\mathcal{M}}(f)}{\longrightarrow} \Upsilon_{\mathcal{M}}(N)\rightarrow F\rightarrow 0, \notag
\end{equation}
for some morphism $f : M\rightarrow N$ in $\mathcal{M}$ that is an admissible morphism in $\mathcal{E}$.
\item[(3)]
Let $\eff(\mathcal{M})$ be the full subcategory of $\Mod (\mathcal{M})$ consisting of all {\em effaceable functors}, i.e. those functors $F$ that admit a projective presentation
\begin{equation}
\Upsilon_{\mathcal{M}}(M)\overset{\Upsilon_{\mathcal{M}}(f)}{\longrightarrow} \Upsilon_{\mathcal{M}}(N)\rightarrow F\rightarrow 0,\notag
\end{equation}
for some morphism $f : M\rightarrow N$ in $\mathcal{M}$ that is a deflation in $\mathcal{E}$.
\item[(4)]
Let $\mathcal{F}(\mathcal{M})$ be the full subcategory of $\Mod (\mathcal{M})$ consisting of those functors $F$ that admit a projective presentation
\begin{equation}
\Upsilon_{\mathcal{M}}(M)\overset{\Upsilon_{\mathcal{M}}(f)}{\longrightarrow} \Upsilon_{\mathcal{M}}(N)\rightarrow F\rightarrow 0,\notag
\end{equation}
for some morphism $f : M\rightarrow N$ in $\mathcal{M}$ that is an admissible morphism in $\mathcal{E}$ and $M\twoheadrightarrow \Imm(f)$ is a right $\mathcal{M}$-approximation.
\item[(5)]
A functor $F\in \Mod(\mathcal{M})$ is called weakly effaceable, if for each object $M\in \mathcal{M}$ and $x\in F(M)$
there exists an admissible epimorphism $f : N \rightarrow M$ such that $F(f)(x) = 0$. We denote by
$\rm Eff(\mathcal{M})$ the full subcategory of all weakly effaceable functors.
\end{itemize}
\end{definition}

\begin{remark}\label{rm2}
Let $\mathcal{M}$ be a contravariantly finite subcategory of an exact category $\mathcal{E}$. A functor $F\in \Mod(\mathcal{M})$ belongs to $\mathcal{F}(\mathcal{M})$, if and only if $F$ admit a projective presentation
\begin{equation}
\Upsilon_{\mathcal{M}}(X)\overset{\Upsilon_{\mathcal{M}}(f)}{\longrightarrow} \Upsilon_{\mathcal{M}}(N)\rightarrow F\rightarrow 0,\notag
\end{equation}
such that $X\in \mathcal{E}$, $N\in \mathcal{M}$ and $f$ is an inflation.
\end{remark}

\begin{lemma}
Let $\mathcal{E}$ be an exact category and $\mathcal{M}$ be an additive subcategory of it. Then $\modd_{\adm}\mathcal{M}$ is an extension-closed subcategory of $\Mod(\mathcal{M})$, and so $\modd_{\adm}\mathcal{M}$ inherits an exact structure from $\Mod(\mathcal{M})$.
\begin{proof}
The proof of \cite[Proposition 3.4]{HRK} carries over.
\end{proof}
\end{lemma}

\begin{proposition}\label{pro2}
Let $\mathcal{E}$ be an exact category and $\mathcal{M}$ be a generating and contravariantly finite subcategory of it. Then $\Upsilon_{\mathcal{M}}$ induces a full and faithful functor
$\Upsilon_{\mathcal{M}}:\mathcal{E}\rightarrow \modd_{\adm}(\mathcal{M})$.
\begin{proof}
For an arbitrary object $X\in \mathcal{E}$, consider conflations $K^0\rightarrowtail M^0\twoheadrightarrow X$ and $K^1\rightarrowtail M^1\twoheadrightarrow K^0$ where $M^0\twoheadrightarrow X$ and $M^1\twoheadrightarrow K^0$ are right $\mathcal{M}$-approximations. Then the induced sequence
\begin{equation}
\Upsilon_{\mathcal{M}}(M^1)\longrightarrow \Upsilon_{\mathcal{M}}(M^0)\longrightarrow \Upsilon_{\mathcal{M}}(X)\rightarrow 0 \notag
\end{equation}
is exact. This implies that $\Upsilon_{\mathcal{M}}(X)\in \modd_{\adm}(\mathcal{M})$. Let $Y\in \mathcal{E}$ and construct a similar projective presentation for $\Upsilon_{\mathcal{M}}(Y)$. By Yoneda's lemma we can see that $\Hom_{\mathcal{E}}(X,Y)\cong \Hom_{\modd_{\adm}(\mathcal{M})}(\Upsilon_{\mathcal{M}}(X),\Upsilon_{\mathcal{M}}(Y))$.
\end{proof}
\end{proposition}

Recall that an exact category $\mathcal{E}$ said to have enough projective objects if for every object $X\in \mathcal{E}$ there exist a projective object $P\in \mathcal{E}$ and a conflation $Y\rightarrowtail P\twoheadrightarrow X$.

In the following proposition we realize exact categories with enough projective objects.

\begin{proposition}\label{pro14}
Let $\mathcal{E}$ be an exact category with enough projective objects $\mathcal{P}$. Then the functor $\mathsf{i}:\mathcal{E}\rightarrow \Mod \mathcal{P}$ given by $\mathsf{i}(X)=\Hom_{\mathcal{E}}(-, X)|_\mathcal{P}$ is full and faithful, it's essential image is $\modd_{adm}\mathcal{P}$, and as exact categories $\mathcal{E}\simeq \modd_{adm}\mathcal{P}$.
\begin{proof}
The proof of the statement that $\mathsf{i}$ is full and faithful is standard. Since $\mathcal{E}$ has enough projectives, $\mathsf{i}(X)\in \modd_{adm}\mathcal{P}$ for every object $X\in \mathcal{E}$. Conversely, let $F\in \modd_{adm}\mathcal{P}$. Choose a projective presentation
\begin{equation}
\Hom_{\mathcal{E}}(-, P)\overset{\Hom_{\mathcal{E}}(-, f)}{\longrightarrow} \Hom_{\mathcal{E}}(-, Q)\rightarrow F\rightarrow 0,\notag
\end{equation}
for $F$ where $f$ is an admissible morphism. Now it is easy to see that $F\cong \Hom_{\mathcal{E}}(-,\Coker(f))|_{\mathcal{P}}$.

For the last statement, first note that $\mathsf{i}$ maps conflations of $\mathcal{E}$ to conflations of $\modd_{adm}\mathcal{P}$. Conversely, let $X\overset{f}{\rightarrow}Y\overset{g}{\rightarrow}Z$ be a complex in $\mathcal{E}$ such that $\mathsf{i}(X)\overset{\mathsf{i}(f)}{\rightarrow}\mathsf{i}(Y)\overset{\mathsf{i}(g)}{\rightarrow}\mathsf{i}(Z)$ is a conflation in $\modd_{adm}\mathcal{P}$. Choose a conflation $K\rightarrowtail P\twoheadrightarrow Z$ with $P\in\mathcal{P}$, and taking the pullback we obtain the following commutative diagram.
\begin{center}
\begin{tikzpicture}
\node (X1) at (0,2) {$K$};
\node (X2) at (2,2) {$K$};
\node (X3) at (-2,0) {$X$};
\node (X4) at (0,0) {$W$};
\node (X5) at (2,0) {$P$};
\node (X6) at (-2,-2) {$X$};
\node (X7) at (0,-2) {$Y$};
\node (X8) at (2,-2) {$Z$};
\draw [-,double,thick] (X1) -- (X2) node [midway,above] {};
\draw [->,thick] (X3) -- (X4) node [midway,above] {};
\draw [->,thick] (X4) -- (X5) node [midway,above] {$w$};
\draw [->,thick] (X6) -- (X7) node [midway,above] {$f$};
\draw [->,thick] (X7) -- (X8) node [midway,above] {$g$};
\draw [>->,thick] (X1) -- (X4) node [midway,left] {};
\draw [>->,thick] (X2) -- (X5) node [midway,left] {};
\draw [-,double,thick] (X3) -- (X6) node [midway,left] {};
\draw [->>,thick] (X4) -- (X7) node [midway,right] {$u$};
\draw [->>,thick] (X5) -- (X8) node [midway,right] {$p$};
\end{tikzpicture}
\end{center}
Since \begin{equation}
\Hom_{\mathcal{E}}(-, X)|_{\mathcal{P}}\overset{\Hom_{\mathcal{E}}(-, f)}{\longrightarrow} \Hom_{\mathcal{E}}(-, Y)|_{\mathcal{P}}\overset{\Hom_{\mathcal{E}}(-, g)}{\longrightarrow} \Hom_{\mathcal{E}}(-, X)|_{\mathcal{P}}\notag
\end{equation}
is a conflation, there exists a morphism $h:P\rightarrow Y$ such that $p1_P=gh$, where $1_P:P\rightarrow P$ is an identity morphism. By the universal property of the pullback, $w$ is a retraction and hence $pw=gu$ is a deflation. Now by the obscure axiom, $g$ is a conflation.
\end{proof}
\end{proposition}

\subsection{Change of categories}
Let $\mathcal{E}$ be an exact category and $\mathcal{M}$ be a full subcategory of $\mathcal{E}$. There is a restriction functor $\mathbb{R}:\Mod (\mathcal{E})\rightarrow\Mod (\mathcal{M})$ that sends any functor $F\in \Mod (\mathcal{E})$ to the composition of $F$ with the inclusion $\mathcal{M}\hookrightarrow \mathcal{E}$.
$\mathbb{R}$ is an exact functor. Also there is a functor
\begin{equation}
\mathbb{E}=\mathcal{E}\otimes_{\mathcal{M}}-:\Mod (\mathcal{M})\longrightarrow\Mod (\mathcal{E})\notag
\end{equation}
defined by $\mathbb{E}(F)(X)=\Hom_{\mathcal{E}}(X,-)|_{\mathcal{M}}\otimes_{\mathcal{M}}F$ for all $F\in \Mod (\mathcal{M})$ and $X\in \mathcal{E}$ \cite[Section 3]{Au1}.

In the following proposition we recall some of the basic properties of $\mathbb{E}$.

\begin{proposition}\label{pro4}
Let $\mathcal{E}$ be an exact category and $\mathcal{M}$ be a full subcategory of $\mathcal{E}$. Then the following statements hold.
\begin{itemize}
\item[(1)]
$\mathbb{E}$ is a right exact functor and preserve sums.
\item[(2)]
The composition $\mathbb{R}\mathbb{E}$ is the identity functor on $\Mod (\mathcal{M})$.
\item[(3)]
For each object $M\in \mathcal{M}$, we have $\mathbb{E}(\Upsilon_{\mathcal{M}}(M))=\Upsilon_{\mathcal{E}}(M)$.
\item[(4)]
$\mathbb{E}$ is a fully faithful functor and a right adjoint for $\mathbb{R}$.
\end{itemize}
\begin{proof}
See \cite[Section 3]{Au1}.
\end{proof}
\end{proposition}

\begin{proposition}\label{pro5}
Let $\mathcal{E}$ be an exact category and $\mathcal{M}$ be an $n$-cluster tilting subcategory of $\mathcal{E}$. Then $\eff(\mathcal{M})=\Eff(\mathcal{M})\cap \modd(\mathcal{M})$.
\begin{proof}
Let $F\in \eff(\mathcal{M})$. Thus there is a projective presentation
$\Upsilon_{\mathcal{M}}(M)\longrightarrow \Upsilon_{\mathcal{M}}(N)\longrightarrow F\rightarrow 0,$
induced by a deflation $f:M\twoheadrightarrow N$ with $M,N\in \mathcal{M}$. Let $X\in \mathcal{M}$ and $x\in F(X)$. $x$ is the image of some morphism $h: X\rightarrow N$. By Proposition 2.9, $f$ sites in an $n$-exact sequence. Now take the $n$-pullback along $h$, we have the following commutative diagram where the horizontal morphisms are deflations.
\begin{center}
\begin{tikzpicture}
\node (X1) at (-2,2) {$X'$};
\node (X2) at (0,2) {$X$};
\node (X3) at (-2,0) {$M$};
\node (X4) at (0,0) {$N$};
\draw [->>,thick] (X1) -- (X2) node [midway,above] {$g$};
\draw [->,thick] (X2) -- (X4) node [midway,right] {h};
\draw [->,thick] (X1) -- (X3) node [midway,left] {};
\draw [->>,thick] (X3) -- (X4) node [midway,below] {f};
\end{tikzpicture}
\end{center}
Obviously $F(g)(x)=0$ and hence $F\in \Eff(\mathcal{M})$.

Now let $F\in \Eff(\mathcal{M})\cap \modd \mathcal{M}$. There is a projective presentation $\Upsilon_{\mathcal{M}}(M)\longrightarrow \Upsilon_{\mathcal{M}}(N)\longrightarrow F\rightarrow 0,$
induced by a morphism $f:M\twoheadrightarrow N$ with $M,N\in \mathcal{M}$. Let $x$ be the image of $1_N$ under the map $\Hom(N,N)\rightarrow F(N)$. By hypothesis there exists a deflation $h:X\twoheadrightarrow N$, such that $F(h)(x) =0$. This means that $h$ must factor through $f$. Thus by Proposition \ref{pro3}, $(0,f):X\oplus M\rightarrow N$ is a deflation and $F\cong \Coker \Upsilon_{\mathcal{M}}(0,f)$.
\end{proof}
\end{proposition}

By Proposition \ref{pro4} and Remark \ref{rm2}, $\mathbb{E}(\eff(\mathcal{M}))\subseteq \eff(\mathcal{E})$ and $\mathbb{E}(\mathcal{F}(\mathcal{M}))\subseteq \mathcal{F}(\mathcal{E})$. In the following proposition we the converse implications under some assumption.

\begin{proposition}\label{pro6}
Let $\mathcal{E}$ be an exact category and $\mathcal{M}$ be an $n$-cluster tilting subcategory of $\mathcal{E}$. Then the following hold.
\begin{itemize}
\item[(1)]
$\mathbb{R}(\modd(\mathcal{E}))\subseteq \modd(\mathcal{M})$.
\item[(2)]
$\mathbb{R}(F)\in \Eff(\mathcal{M})$ if and only if $F\in \Eff(\mathcal{E})$.
\item[(3)]
$\mathbb{R}(\eff(\mathcal{E}))\subseteq \eff(\mathcal{M})$.
\end{itemize}
\begin{proof}
Let $F\in \modd(\mathcal{E})$ and consider a projective presentation
\begin{equation}
\Upsilon_{\mathcal{E}}(X)\overset{\Upsilon_{\mathcal{E}}(f)}{\longrightarrow} \Upsilon_{\mathcal{E}}(Y)\rightarrow F\rightarrow 0,\notag
\end{equation}
such that $f:X\rightarrow Y$ is a deflation in $\mathcal{E}$. Because $\mathcal{M}$ is generating and contravariantly finite, there exist acyclic complexes $M^1\rightarrow M^0\rightarrow X\rightarrow 0$ and $N^1\rightarrow N^0\rightarrow Y\rightarrow 0$ with $M^i,N^i\in \mathcal{M}$ for each $i=1,2$, such that we have the following exact sequences.
\begin{align*}
&\Upsilon_{\mathcal{M}}(M^1)\rightarrow\Upsilon_{\mathcal{M}}(M^0)\rightarrow\Upsilon_{\mathcal{M}}(X)\rightarrow 0\\
&\Upsilon_{\mathcal{M}}(N^1)\rightarrow\Upsilon_{\mathcal{M}}(N^0)\rightarrow\Upsilon_{\mathcal{M}}(Y)\rightarrow 0.
\end{align*}
Now because $\mathbb{R}$ is an exact functor, we obtain the following commutative diagram with exact columns and exact row.
\begin{center}
\begin{tikzpicture}
\node (X1) at (-3,4) {$\Upsilon_{\mathcal{M}}(M^1)$};
\node (X2) at (0,4) {$\Upsilon_{\mathcal{M}}(N^1)$};
\node (X3) at (-3,2) {$\Upsilon_{\mathcal{M}}(M^0)$};
\node (X4) at (0,2) {$\Upsilon_{\mathcal{M}}(N^0)$};
\node (X5) at (-3,0) {$\Upsilon_{\mathcal{M}}(X)$};
\node (X6) at (0,0) {$\Upsilon_{\mathcal{M}}(Y)$};
\node (X7) at (2,0) {$F|_{\mathcal{M}}$};
\node (X8) at (4,0) {$0$};
\node (X9) at (-3,-2) {$0$};
\node (X10) at (0,-2) {$0$};
\draw [->,thick] (X1) -- (X2) node [midway,above] {};
\draw [->,thick] (X3) -- (X4) node [midway,above] {};
\draw [->,thick] (X5) -- (X6) node [midway,above] {};
\draw [->,thick] (X6) -- (X7) node [midway,above] {};
\draw [->,thick] (X7) -- (X8) node [midway,above] {};
\draw [->,thick] (X1) -- (X3) node [midway,left] {};
\draw [->,thick] (X3) -- (X5) node [midway,left] {};
\draw [->,thick] (X5) -- (X9) node [midway,left] {};
\draw [->,thick] (X2) -- (X4) node [midway,left] {};
\draw [->,thick] (X4) -- (X6) node [midway,left] {};
\draw [->,thick] (X6) -- (X10) node [midway,left] {};
\end{tikzpicture}
\end{center}
Now it is easy to check that the induced sequence
\begin{equation}
\Upsilon_{\mathcal{M}}(N^1)\oplus \Upsilon_{\mathcal{M}}(M^0)\longrightarrow \Upsilon_{\mathcal{M}}(N^0) \longrightarrow F|_{\mathcal{M}}\rightarrow 0\notag
\end{equation}
is exact. This proves $(1)$.

Now let $F\in \Eff(\mathcal{E})$, and denote by $\widetilde{F}$ the restriction functor $F|_{\mathcal{M}}=\mathbb{R}(F)$. Let $M\in \mathcal{M}$ and $x\in \widetilde{F}(M)=F(M)$. By assumption there exists a deflation $f:E\twoheadrightarrow M$ such that $F(f)(x)=0$. Since $\mathcal{M}$ is generating, there is a deflation $g:N\twoheadrightarrow E$ for some $N\in \mathcal{M}$. Now obviously $\widetilde{F}(fg)(x)=0$. Conversely, assume that $\mathbb{R}(F)\in \Eff(\mathcal{M})$, $X\in \mathcal{M}$ and $x\in F(X)$. Because $\mathcal{M}$ is a generating subcategory, there is a deflation $f:M\twoheadrightarrow X$ for some $M\in \mathcal{M}$. The assumption that $\widetilde{F}$ is weakly effaceable implies that there is a deflation $g:N\twoheadrightarrow M$ such that $\widetilde{F}(g)\big(F(f)(x)\big)=0$, which proves $(2)$.

$(3)$ follows from $(1)$, $(2)$ and Proposition \ref{pro5}.
\end{proof}
\end{proposition}

\begin{proposition}
Let $\mathcal{E}$ be an exact category and $\mathcal{M}$ be an $n$-cluster tilting subcategory of $\mathcal{E}$.
\begin{itemize}
\item[(1)]
$\Eff(\mathcal{M})$ is a localising subcategory of $\Mod(\mathcal{M})$
\item[(2)]
The restriction functor
$\mathbb{R}:\Mod(\mathcal{E})\rightarrow \Mod(\mathcal{M})$ induces an equivalence $\dfrac{\Mod(\mathcal{E})}{\Eff(\mathcal{E})}\overset{\widehat{\mathbb{R}}}{\simeq}\dfrac{\Mod(\mathcal{M})}{\Eff(\mathcal{M})}$ making the following diagram commutative.
\begin{center}
\begin{tikzpicture}
\node (X1) at (-3,3) {$\Mod(\mathcal{E})$};
\node (X2) at (0,3) {$\Mod(\mathcal{M})$};
\node (X3) at (-3,0) {$\dfrac{\Mod(\mathcal{E})}{\Eff(\mathcal{E})}$};
\node (X4) at (0,0) {$\dfrac{\Mod(\mathcal{M})}{\Eff(\mathcal{M})}$};
\draw [->,thick] (X1) -- (X2) node [midway,above] {$\mathbb{R}$};
\draw [->,thick] (X2) -- (X4) node [midway,right] {$\mathbb{Q}_{\mathcal{M}}$};
\draw [->,thick] (X1) -- (X3) node [midway,left] {$\mathbb{Q}_{\mathcal{E}}$};
\draw [->,thick] (X3) -- (X4) node [midway,below] {$\widehat{\mathbb{R}}$};
\end{tikzpicture}
\end{center}
\end{itemize}
\begin{proof}
The first statement follows from \cite[Proposition 3.5]{E}. Thus $\mathbb{Q}_{\mathcal{M}}$ has a fully faithful right adjoint, so the composition
\begin{equation}
\Mod(\mathcal{E})\overset{\mathbb{R}}{\longrightarrow} \Mod(\mathcal{M})\overset{\mathbb{Q}_{\mathcal{M}}}{\longrightarrow}\dfrac{\Mod(\mathcal{M})}{\Eff(\mathcal{M})} \notag
\end{equation}
is an exact functor and has a fully faithful right adjoint and by Proposition \ref{pro6}. The kernel of this exact functor is equal to $\Eff(\mathcal{E})$, then the result follow from \cite{G}.
\end{proof}
\end{proposition}

\subsection{Homological properties of $\modd_{\adm}(\mathcal{M})$}
Let $\mathcal{E}$ be an exact category and $\mathcal{M}$ be an $n$-cluster tilting subcategory of $\mathcal{E}$. In this subsection we study the exact category $\modd_{\adm}(\mathcal{M})$ and show that $\modd_{\adm}(\mathcal{M})$ is an $n$-Auslander exact category which is the higher dimensional analogue of Auslander exact category introduced in \cite[Section 4.1]{HRK}.
\begin{lemma}
Let $E\in \eff(\mathcal{M})$ and let $\Upsilon_{\mathcal{M}}(M)\overset{\Upsilon_{\mathcal{M}}(f)}{\longrightarrow}\Upsilon_{\mathcal{M}}(N)\rightarrow E\rightarrow 0$ be a projective presentation of $E$. Then there is an $M'\in \mathcal{M}$ such that $(0,f):M'\oplus M\twoheadrightarrow N$ is a deflation.
\begin{proof}
By Proposition \ref{pro5}, $E\in \Eff(\mathcal{M})$. Assume that $x$ is the image of $1_N$ under the map $\Hom_{\mathcal{E}}(N,N)|_{\mathcal{M}}\rightarrow E(N)$. By hypothesis there exists a deflation $h:M'\twoheadrightarrow N$ such that $F(h)(x) =0$. This means that $h$ must factor through $f$. Thus by Proposition \ref{pro3}, $(0,f):M'\oplus M\rightarrow N$ is a deflation and $F\cong \Coker(\Upsilon_{\mathcal{M}}(0,f))$.
\end{proof}
\end{lemma}

\begin{proposition}\label{pro7}
The category $\modd_{\adm}(\mathcal{M})$ is obtained as $\mathbb{Q}_{\mathcal{M}}^{-1}(\Upsilon_{\mathcal{M}}(\mathcal{E}))\cap \modd\mathcal{M}$.
\begin{proof}
Let $F\in \mathbb{Q}_{\mathcal{M}}^{-1}(\Upsilon_{\mathcal{M}}(\mathcal{E}))\cap \modd\mathcal{M}$. There is a projective presentation
\begin{center}
$\Upsilon_{\mathcal{M}}(A)\overset{\Upsilon_{\mathcal{M}}(f)}{\longrightarrow}\Upsilon_{\mathcal{M}}(B)\rightarrow F\rightarrow 0.$
\end{center}
Assume that $\mathbb{Q}_{\mathcal{M}}(F)\cong \Upsilon_{\mathcal{M}}(C)$ for some object $C\in \mathcal{E}$. This means that there is a morphism $\eta:F\rightarrow \Upsilon_{\mathcal{M}}(C)$ such that $\Ker(\eta),\Coker(\eta)\in \Eff(\mathcal{M})$. By the Yoneda's lemma the composition $\Upsilon_{\mathcal{M}}(B)\rightarrow F\rightarrow \Upsilon_{\mathcal{M}}(C)$ is induced by a morphism $g:B\rightarrow C$ and it's cokernel is equivalent to $\Coker(\eta)$ and so is effaceable.
The similar argument as in the proof of \cite[Proposition 3.23]{HRK} shows that, there is an object $B'$ such that
\begin{equation}
\widetilde{f}=\begin{pmatrix}
1_{B'}&0\\
0&f
\end{pmatrix}:B'\oplus A\longrightarrow B'\oplus B \notag
\end{equation}
is admissible and obviously $\Coker(\Upsilon_{\mathcal{M}}(\widetilde{f}))\cong F$. Since $\mathcal{M}$ is a generating subcategory, $B'$ can be chosen from $\mathcal{M}$ and so $F\in \modd_{\adm}(\mathcal{M})$.
Now let $F\in \modd_{\adm}(\mathcal{M})$. Choose a projective presentation
\begin{equation}
\Hom_{\mathcal{E}}(-, X)|_{\mathcal{M}}\overset{\Hom_{\mathcal{E}}(-, f)}{\longrightarrow} \Hom_{\mathcal{E}}(-, Y)|_{\mathcal{M}}\rightarrow F\rightarrow 0,\notag
\end{equation}
for $F$ where $f$ is an admissible morphism in $\mathcal{E}$. Then $F\cong \Hom_{\mathcal{E}}(-, \Coker(f))|_{\mathcal{M}}$ and the result follows.
\end{proof}
\end{proposition}

\begin{lemma}\label{lem1}
Let $\mathcal{E}$ be an exact category and $\mathcal{M}$ be an $n$-cluster tilting subcategory of $\mathcal{E}$. Then
the adjoint pair $(\mathbb{R},\mathbb{E}):\Mod(\mathcal{E})\rightarrow \Mod(\mathcal{M})$ can be restricted to an adjoint pair $(\mathbb{R},\mathbb{E}):\modd_{\adm}(\mathcal{E})\rightarrow \modd_{\adm}(\mathcal{M})$.
\begin{proof}
Obviously $\mathbb{E}(\modd_{\adm}(\mathcal{M}))\subseteq \modd_{\adm}(\mathcal{E})$. Now let $F\in \modd_{\adm}(\mathcal{E})$. By \cite[Proposition 3.23]{HRK}, $\mathbb{Q}_{\mathcal{E}}(F)\simeq \Upsilon_{\mathcal{E}}(C)$ in $\dfrac{\Mod(\mathcal{E})}{\Eff(\mathcal{E})}$ for some $C\in \mathcal{E}$. Thus  $\mathbb{Q}_{\mathcal{M}}(\mathbb{R}(F))\simeq \Upsilon_{\mathcal{M}}(C)$, and because $\mathbb{R}(F)\in \modd \mathcal{M}$, by Proposition \ref{pro7} we have that $\mathbb{R}(F)\in \modd_{\adm}\mathcal{M}$.
\end{proof}
\end{lemma}

\begin{remark}
Let $\mathcal{E}$ be an exact category and $\mathcal{M}$ be a generating and contravariantly finite subcategory of it. Then the induced functor
$\Upsilon_{\mathcal{M}}:\mathcal{E}\rightarrow \modd_{\adm}(\mathcal{M})$ (see Proposition \ref{pro2}) is a left exact functor.
For a conflation $X\rightarrowtail Y\twoheadrightarrow Z$ in $\mathcal{M}$ we have the exact sequence
\begin{equation}
0\rightarrow \Upsilon_{\mathcal{M}}(X)\longrightarrow \Upsilon_{\mathcal{M}}(Y)\longrightarrow \Upsilon_{\mathcal{M}}(Z) \notag
\end{equation}
in $\Mod(\mathcal{M})$. Lemma \ref{lem1} implies that the image of $ \Upsilon_{\mathcal{M}}(Y)\longrightarrow \Upsilon_{\mathcal{M}}(Z)$ belongs to $\modd_{\adm}(\mathcal{M})$ which shows that the functor $\Upsilon_{\mathcal{M}}:\mathcal{E}\rightarrow \modd_{\adm}(\mathcal{M})$ is left exact.
\end{remark}

\begin{proposition}\label{pro8}
$(\eff(\mathcal{M}),\mathcal{F}(\mathcal{M}))$ is a torsion pair in $\modd_{\adm}\mathcal{M}$ such that $\eff(\mathcal{M})$ satisfies axiom $\textbf{A2}$.
\begin{proof}
We prove the result using Proposition \ref{pro7} and the similar result about the torsion pair $(\eff(\mathcal{E}),\mathcal{F}(\mathcal{E}))$ in $\modd_{\adm}\mathcal{E}$ \cite[Proposition 3.5]{HRK}.

By Proposition \ref{pro4}, $\mathbb{E}$ is full and faithful. By Proposition \ref{pro4} and Remark \ref{rm2}, $\mathbb{E}(\eff(\mathcal{M}))\\\subseteq \eff(\mathcal{E})$ and $\mathbb{E}(\mathcal{F}(\mathcal{M}))\subseteq \mathcal{F}(\mathcal{E})$. Then by \cite[Proposition 3.5]{HRK}, we have that\\ $\Hom_{\modd_{\adm}\mathcal{M}}(\eff(\mathcal{M}),\mathcal{F}(\mathcal{M}))=0$. Now let $G\in \modd_{\adm}\mathcal{M}$ and choose a projective presentation $\Upsilon(M)\overset{\Upsilon(f)}{\rightarrow} \Upsilon(N)\rightarrow G\rightarrow 0$ where $f$ is an admissible morphism. By the proof of \cite[Proposition 3.5]{HRK}, there is an exact sequence $0\rightarrow F\rightarrow \mathbb{E}(G)\rightarrow H\rightarrow 0$ in $\Mod\mathcal{E}$ where $F$ and $H$ have projective presentations
\begin{align*}
&\Hom_{\mathcal{E}}(-,M)\longrightarrow \Hom_{\mathcal{E}}(-,\Imm(f))\longrightarrow F\longrightarrow 0,\\
&\Hom_{\mathcal{E}}(-,\Imm(f))\longrightarrow \Hom_{\mathcal{E}}(-,N)\longrightarrow H\longrightarrow 0.
\end{align*}
By applying the functor $\mathbb{R}$ we see that $0\rightarrow \mathbb{R}(F)\rightarrow G\rightarrow \mathbb{R}(H)\rightarrow 0$ is an exact sequence with $\mathbb{R}(F)\in \eff(\mathcal{M})$ and $\mathbb{R}(H)\in \mathcal{F}(\mathcal{M})$. Therefore $(\eff(\mathcal{M}),\mathcal{F}(\mathcal{M}))$ is a torsion pair in $\modd_{\adm}\mathcal{M}$. Let $F\in \modd_{\adm}\mathcal{M}$, $G\in \eff(\mathcal{M})$ and $\eta:F\rightarrow G$ be an arbitrary morphism. Obviously $\mathbb{E}(F)\in \modd_{\adm}\mathcal{E}$ and $\mathbb{E}(G)\in \eff(\mathcal{E})$. By \cite[Proposition 3.5]{HRK}, Lemma \ref{lem1} and Proposition \ref{pro6}, $\eta$ is admissible with image in $\eff(\mathcal{M})$ and the result follows.
\end{proof}
\end{proposition}

\begin{proposition}\label{pro9}
Let $\mathcal{E}$ be an exact category and $\mathcal{M}$ be a generating contravariantly finite subcategory of $\mathcal{E}$. Then $\Upsilon_{\mathcal{M}}:\mathcal{E}\rightarrow \modd_{\adm}(\mathcal{M})$ admit a left adjoint $L_{\mathcal{M}}:\modd_{\adm}(\mathcal{M})\\\rightarrow \mathcal{E}$. Moreover $L_{\mathcal{M}}$ is an exact functor and induces an equivalence of exact categories $\dfrac{\modd_{\adm}(\mathcal{M})}{\eff(\mathcal{M})}\simeq \mathcal{E}$.
\begin{proof}
The proof is similar to the proof of \cite[Theorem 3.7]{HRK}. For each $F\in \modd_{\adm}(\mathcal{M})$ we choose a projective presentation $\Upsilon_{\mathcal{M}}(M)\overset{\Upsilon_{\mathcal{M}}(f)}{\longrightarrow}\Upsilon_{\mathcal{M}}(N)\rightarrow F\rightarrow 0$. If $F=\Upsilon_{\mathcal{M}}(M)$ we simply choose $0\rightarrow \Upsilon_{\mathcal{M}}(M)\rightarrow \Upsilon_{\mathcal{M}}(M)\rightarrow 0$. Then we have $L_{\mathcal{M}}(F)=\Coker(f)$. Obviously $\Ker(L_{\mathcal{M}})=\eff(\mathcal{M})$, thus by \cite[Proposition 1.3]{GZ} we have desired equivalence.
\end{proof}
\end{proposition}

\begin{definition}$($\cite[Definition 3.13]{HRK}$)$
Let $\mathcal{E}$ be an exact category and $\mathcal{X}$ be a full subcategory of $\mathcal{E}$. $\cogen(\mathcal{X})$ is the full subcategory of $\mathcal{E}$ consisting of those objects $E$ such that there is an inflation $E\rightarrowtail X$ for some object $X\in\mathcal{X}$. $\gen(\mathcal{X})$ is defined dually.
\end{definition}

Let $\mathcal{A}$ be an additive category and $\mathcal{X}$ be a class of objects in $\mathcal{A}$. In the
following we write $^{\bot}\mathcal{X}:=\{A\in\mathcal{A}|\Hom_{\mathcal{A}}(A, X)=0, \forall X\in \mathcal{X}\}$.

The following proposition is a higher dimensional version of \cite[Proposition 3.14]{HRK}

\begin{proposition}\label{pro12}
Consider the torsion pair $(\eff(\mathcal{M}),\mathcal{F}(\mathcal{M}))$ in $\modd_{\adm}(\mathcal{M})$ as in Proposition \ref{pro8} and set $\mathcal{Q}=\proj(\modd_{\adm}(\mathcal{M}))$. Then the following hold.
\begin{itemize}
\item[(1)]
$\eff(\mathcal{M})=^{\bot}\mathcal{Q}$.
\item[(2)]
The Yoneda functor $H:\mathcal{M}\rightarrow \modd_{\adm}\mathcal{M}$ gives an equivalence $\mathcal{M}\simeq \mathcal{Q}$.
\item[(3)]
$\mathcal{F}(\mathcal{M})=\cogen(\mathcal{Q})$.
\item[(4)]
$\Ext_{\modd_{\adm}\mathcal{M}}^{1,\cdots,n}(\eff(\mathcal{M}),\mathcal{Q})=0$.
\end{itemize}
\begin{proof}
$(1).$ It follows from \cite[Proposition 3.14]{HRK} and the properties of the adjoint pair $(\mathbb{R},\mathbb{E})$.

$(2).$ It follows from $(1)$.

$(3).$ The proof is similar to the proof of \cite[Proposition 3.14]{HRK}.

$(4).$ let $E\in \eff(\mathcal{M})$. Then there is a deflation $M^{n}\twoheadrightarrow M^{n+1}$ that induces a projective presentation of $E$. Since $\mathcal{M}$ is an $n$-cluster tilting subcategory of $\mathcal{E}$, by \cite[Theorem 4.14]{J}, there is an $n$-exact sequence
\begin{equation}\label{ne}
M^0\rightarrowtail M^1\rightarrow \cdots\rightarrow M^n\twoheadrightarrow M^n,
\end{equation}
in $\mathcal{M}$. Since $\mathcal{M}$ is $n$-rigid,
\begin{equation}
0\rightarrow \Upsilon_{\mathcal{M}}(M^0)\rightarrow\Upsilon_{\mathcal{M}}(M^1)\rightarrow\cdots\Upsilon_{\mathcal{M}}(M^n)\rightarrow\Upsilon_{\mathcal{M}}(M^{n+1})\rightarrow E\rightarrow 0, \notag
\end{equation}
is a projective resolution of $E$. Now let $\Upsilon_{\mathcal{M}}(X)\in \proj(\modd_{\adm}(\mathcal{M}))$. Applying the functor $\Hom_{\modd_{\adm}(\mathcal{M})}(-, \Upsilon_{\mathcal{M}}(X))$ to this projective resolution. By using Yoneda's lemma and the fact that \eqref{ne} is an $n$-exact sequence the result follows.
\end{proof}
\end{proposition}

\begin{definition}$($\cite[Definition 4.15]{HRK}$)$
Let $\mathcal{X}$ be a subcategory of an exact category $\mathcal{E}$. A morphism $f:E\rightarrow X$ in $\mathcal{E}$ with
$X\in \mathcal{X}$ is called an {\em admissible left $\mathcal{X}$-approximation} if $f$ is an admissible morphism and any morphism
$E\rightarrow X'$ with $X'\in \mathcal{X}$ factors through $f$. The subcategory $\mathcal{X}$ is called {\em admissibly covariantly finite} if for
all objects $E$ in $\mathcal{E}$ there exists an admissible left $\mathcal{X}$-approximation $f:E\rightarrow X$. The notions
of {\em admissibly contravariantly finite subcategory} of $\mathcal{E}$ and {\em admissible right $\mathcal{X}$-approximation} are defined dually.
\end{definition}

\begin{proposition}\label{pro13}
In the setting of Proposition \ref{pro9}, the following statements hold.
\begin{itemize}
\item[(1)]
For every $F\in \modd_{\adm}(\mathcal{M})$ the unit of adjunction $F\rightarrow \Upsilon_{\mathcal{M}}(L_{\mathcal{M}}(F))$ is an admissible morphism.
\item[(2)]
$\mathcal{Q}=\proj(\modd_{\adm}(\mathcal{M}))$ is an admissibly covariantly finite subcategory of $\modd_{\adm}(\mathcal{M})$.
\end{itemize}
\begin{proof}
$(1).$ Choose a projective presentation
\begin{center}
$\Upsilon_{\mathcal{M}}(M)\overset{\Upsilon_{\mathcal{M}}(f)}{\longrightarrow}\Upsilon_{\mathcal{M}}(N)\rightarrow F\rightarrow 0$
\end{center}
for $F$ and set $C:=\Coker(f)$. We can easily see that the unit is the unique morphism $\eta:F\rightarrow \Upsilon_{\mathcal{M}}(C)$ that makes the following diagram commutative.
\begin{center}
\begin{tikzpicture}
\node (X1) at (-3,0) {$\Upsilon_{\mathcal{M}}(M)$};
\node (X2) at (0,0) {$\Upsilon_{\mathcal{M}}(N)$};
\node (X3) at (2,0) {$F$};
\node (X4) at (4,0) {$0.$};
\node (X5) at (0,-2) {$\Upsilon_{\mathcal{M}}(C)$};
\draw [->,thick] (X1) -- (X2) node [midway,above] {$\Upsilon_{\mathcal{M}}(f)$};
\draw [->,thick] (X2) -- (X3) node [midway,right] {};
\draw [->,thick] (X3) -- (X4) node [midway,left] {};
\draw [->,thick] (X2) -- (X5) node [midway,above] {};
\draw [->,thick,dashed] (X3) -- (X5) node [midway,below] {$\eta$};
\end{tikzpicture}
\end{center}
Applying $L_{\mathcal{M}}$ we obtain the following commutative diagram with exact rows.
\begin{center}
\begin{tikzpicture}
\node (X1) at (-2,0) {$M$};
\node (X2) at (0,0) {$N$};
\node (X3) at (2,0) {$L_{\mathcal{M}}(F)=C$};
\node (X4) at (4,0) {$0$};
\node (X5) at (-2,-2) {$M$};
\node (X6) at (0,-2) {$N$};
\node (X7) at (2,-2) {$C$};
\node (X8) at (4,-2) {$0.$};
\draw [->,thick] (X1) -- (X2) node [midway,above] {};
\draw [->,thick] (X2) -- (X3) node [midway,right] {};
\draw [->,thick] (X3) -- (X4) node [midway,left] {};
\draw [->,thick] (X5) -- (X6) node [midway,above] {};
\draw [->,thick] (X6) -- (X7) node [midway,right] {};
\draw [->,thick] (X7) -- (X8) node [midway,left] {};
\draw [-,thick,double] (X1) -- (X5) node [midway,above] {};
\draw [-,thick,double] (X2) -- (X6) node [midway,right] {};
\draw [->,thick] (X3) -- (X7) node [midway,right] {$L_{\mathcal{M}}(\eta)$};
\end{tikzpicture}
\end{center}
Thus $L_{\mathcal{M}}(\eta)$ is an isomorphism. Since any localisation functor reflect admissible morphisms, $\eta$ is an admissible morphism and the result follows.

$(2).$ Let $F\in \modd_{\adm}(\mathcal{M})$. Choose an inflation left $\mathcal{M}$-approximation $L_{\mathcal{M}}(F)\rightarrowtail M$ for $L_\mathcal{M}(F)$. Because $\Upsilon_{\mathcal{M}}$ is a left exact functor $\Upsilon_{\mathcal{M}}(L_\mathcal{M}(F))\rightarrowtail \Upsilon_{\mathcal{M}}(M)$ is also an inflation.
Then obviously composition with unit of adjunction $F\overset{\eta}{\rightarrow}\Upsilon_{\mathcal{M}}(L_{\mathcal{M}}(F))\rightarrowtail \Upsilon_{\mathcal{M}}(M)$ is an admissible left $\mathcal{Q}$-approximation of $F$.

\end{proof}
\end{proposition}

\section{Higher Auslander correspondence}
In this section we define $n$-Auslander exact categories and show that these categories are in correspondence with $n$-cluster tilting subcategories of exact categories up to equivalence of categories.

Let $\mathcal{E}$ be an exact category, $k$ a positive integer and $A, C\in \mathcal{A}$. An exact sequence
\begin{equation}
0\rightarrow A\rightarrow X_{k-1}\rightarrow \cdots \rightarrow X_0\rightarrow C\rightarrow 0 \notag
\end{equation}
in $\mathcal{A}$ is called a {\em $k$-fold extension of $C$ by $A$}. Two $k$-fold extensions of $C$ by $A$,
\begin{equation}
\xi:0\rightarrow A\rightarrow B_{k-1}\rightarrow \cdots \rightarrow B_0\rightarrow C\rightarrow 0 \notag
\end{equation}
and
\begin{equation}
\xi':0\rightarrow A\rightarrow B'_{k-1}\rightarrow \cdots \rightarrow B'_0\rightarrow C\rightarrow 0 \notag
\end{equation}
are said to be {\em Yoneda equivalent} if there is a chain of $k$-fold extensions of $C$ by $A$
\begin{equation}
\xi=\xi_0,  \xi_1, \ldots ,\xi_{l-1}, \xi_l=\xi' \notag
\end{equation}
such that for every $i\in \{0,\cdots,l-1\}$, we have either a chain map $\xi_i\rightarrow \xi_{i+1}$ starting with $1_A$ and ending with $1_C$, or a chain map $\xi_{i+1}\rightarrow \xi_{i}$ starting with $1_A$ and ending with $1_C$.
$\Ext_{\mathcal{A}}^k(C,A)$ is defined as the set of Yoneda equivalence classes of $k$-fold extensions of $C$ by $A$ \cite[Definition 6.1]{FD} (see also \cite{Mac, Mi}).

The following definition is higher dimensional analogue of \cite[Definition 4.1]{HRK}.

\begin{definition}\label{def}
Let $\mathcal{E}$ be an exact category with enough projectives $\mathcal{P}=\proj(\mathcal{E})$ and $n$ be a positive integer. We say that $\mathcal{E}$ is an {\em $n$-Auslander exact category} if it satisfies the following conditions:
\begin{itemize}
\item[(1)]
$(^{\bot}\mathcal{P},\cogen(\mathcal{P}))$ is a torsion pair in $\mathcal{E}$.
\item[(2)]
$^{\bot}\mathcal{P}$ satisfies axiom $\textbf{A2}$.
\item[(3)]
$\Ext_{\mathcal{E}}^{1,\cdots,n}(^{\bot}\mathcal{P},\mathcal{P})=0$.
\item[(4)]
$\gl.dim(\mathcal{E})\leq n+1$.
\item[(5)]
$\mathcal{P}$ is an admissibly covariantly finite subcategory of $\mathcal{E}$.
\end{itemize}
\end{definition}

Note that, by \cite[Lemma 4.16]{HRK}, $1$-Auslander exact categories are precisely Auslander exact categories defined in \cite[Definition 4.1]{HRK}.

\begin{lemma}\label{lem2}
Let $\mathcal{E}$ be an $n$-Auslander exact category with $\mathcal{P}=\proj(\mathcal{E})$. For every $M\in \cogen(\mathcal{P})$, there exists a conflation $M\rightarrowtail \overline{M}\twoheadrightarrow T$ such that $\pd(\overline{M})\leq n-1$ and $T\in ^{\bot}\mathcal{P}$.
\begin{proof}
The proof is similar to the proof of \cite[Proposition 4.4 (3)]{HRK}.
\end{proof}
\end{lemma}

\begin{proposition}\label{pro10}
Let $\mathcal{E}$ be an $n$-Auslander exact category with $\mathcal{P}=\proj(\mathcal{E})$. The composition $H:\mathcal{P}\hookrightarrow \mathcal{E}\rightarrow \frac{\mathcal{E}}{^{\bot}\mathcal{P}}$ is full and faithful, and it's essential image is an $n$-rigid subcategory.
\begin{proof}
The prove of fullness and faithfulness is similar to \cite[Proposition 4.5]{HRK}. We show that the essential image is $n$-rigid. The proof is an adaptation of the proof of \cite[Theorem 4.6]{EN}.

For simplicity we set $\mathcal{E}':=\frac{\mathcal{E}}{^{\bot}\mathcal{P}}$ and $\mathcal{M}:=H(\mathcal{P})$.
Let $k\in \{1,\ldots,n-1\}$, $M_1, M_2\in \mathcal{M}$ and
\begin{center}
$\xi: 0\rightarrow M_2\rightarrow F^{n-k+1}\rightarrow \cdots \rightarrow F^n\rightarrow M_1\rightarrow 0$
\end{center}
be a $k$-fold extension of $M_1$ by $M_2$ in $\mathcal{E'}$. There exists $X',X^{n+1}\in \mathcal{M}$ such that $H(X^{n+1})=M_1$ and $H(X')=M_2$.
Then the cokernel of $F^n\rightarrow H(X^{n+1})$ formed in $\mathcal{E}$, denoted by $C$, belongs to $^{\bot}\mathcal{P}$.
Since $\mathcal{E}$ has enough projectives, there is an epimorphism $H(X^n)\rightarrow F^n$ in $\mathcal{E}'$ for some $X^n\in \mathcal{P}$. Obviously the cokernel of the composition $H(X^n)\rightarrow F^n\rightarrow H(X^{n+1})$ is equal to $C$.
Because $\gl.dim\mathcal{E}\leq n+1$, we can construct projective resolution
\begin{equation}
0\rightarrow X^0\rightarrow X^1\rightarrow \cdots \rightarrow X^n\rightarrow X^{n+1}\rightarrow C\rightarrow 0 \notag
\end{equation}
for $C$. Since $\Ext_{E}^{1,\cdots,n}(^{\bot}\mathcal{P},\mathcal{P})=0$, this projective resolution remains exact after applying $\Hom(-,X)$ for each $X\in \mathcal{P}$. This means that
$X^0\rightarrow X^1\rightarrow \cdots \rightarrow X^n\rightarrow X^{n+1}$ is an $n$-exact sequence in $\mathcal{P}$.
Then we have the following commutative diagram with exact rows in $\mathcal{E}'$.
\begin{center}
\begin{tikzpicture}
\node (X1) at (-10,1) {$0$};
\node (X2) at (-8.5,1) {$H(X^0)$};
\node (X3) at (-7,1) {$\cdots$};
\node (X4) at (-5,1) {$H(X^{n-k})$};
\node (X5) at (-2,1) {$H(X^{n-k+1})$};
\node (X6) at (0,1) {$\cdots$};
\node (X7) at (2,1) {$H(X^n)$};
\node (X8) at (4,1) {$H(X^{n+1})$};
\node (X9) at (5.75,1) {$0$};
\node (X10) at (-6,-0.5) {$0$};
\node (X11) at (-4,-0.5) {$H(X')$};
\node (X12) at (-2,-0.5) {$F^{n-k+1}$};
\node (X13) at (0,-0.5) {$\cdots$};
\node (X14) at (2,-0.5) {$F^n$};
\node (X15) at (4,-0.5) {$H(X^{n+1})$};
\node (X16) at (5.75,-0.5) {$0.$};
\draw [->,thick] (X1) -- (X2) node [midway,left] {};
\draw [->,thick] (X2) -- (X3) node [midway,left] {};
\draw [->,thick] (X3) -- (X4) node [midway,left] {};
\draw [->,thick] (X4) -- (X5) node [midway,left] {};
\draw [->,thick] (X5) -- (X6) node [midway,left] {};
\draw [->,thick] (X6) -- (X7) node [midway,above] {};
\draw [->,thick] (X7) -- (X8) node [midway,left] {};
\draw [->,thick] (X8) -- (X9) node [midway,above] {};
\draw [->,thick] (X10) -- (X11) node [midway,left] {};
\draw [->,thick] (X11) -- (X12) node [midway,above] {};
\draw [->,thick] (X12) -- (X13) node [midway,above] {};
\draw [->,thick] (X13) -- (X14) node [midway,above] {};
\draw [->,thick] (X14) -- (X15) node [midway,above] {};
\draw [->,thick] (X15) -- (X16) node [midway,above] {};
\draw [->,thick] (X7) -- (X14) node [midway,above] {};
\draw [double,-,thick] (X8) -- (X15) node [midway,above] {};
\end{tikzpicture}
\end{center}

The rest of the proof is exactly similar to the proof of \cite[Theorem 4.6]{EN}, for the convenience of the reader we complete the proof.

First assume that $k=1$. By the universal property of kernel there exists a morphism $H(X^{n-1})\rightarrow H(X')$ that makes the solid part of the following diagram commutative.
\begin{center}
\begin{tikzpicture}
\node (X3) at (-9,0.5) {$0$};
\node (X4) at (-7,0.5) {$H(X^0)$};
\node (X5) at (-5.5,0.5) {$\cdots$};
\node (X6) at (-3.5,0.5) {$H(X^{n-1})$};
\node (X7) at (-1.5,0.5) {$H(X^n)$};
\node (X8) at (0.5,0.5) {$H(X^{n+1})$};
\node (X9) at (2.5,0.5) {$0$};
\node (X10) at (-5,-1) {$0$};
\node (X11) at (-3.5,-1) {$H(X')$};
\node (X12) at (-1.5,-1) {$F^n$};
\node (X13) at (0.5,-1) {$H(X^{n+1})$};
\node (X14) at (2.5,-1) {$0.$};
\draw [->,thick] (X3) -- (X4) node [midway,left] {};
\draw [->,thick] (X4) -- (X5) node [midway,left] {};
\draw [->,thick] (X5) -- (X6) node [midway,left] {};
\draw [->,thick] (X6) -- (X7) node [midway,above] {};
\draw [->,thick] (X7) -- (X8) node [midway,left] {};
\draw [->,thick] (X8) -- (X9) node [midway,above] {};
\draw [->,thick] (X10) -- (X11) node [midway,left] {};
\draw [->,thick] (X11) -- (X12) node [midway,above] {};
\draw [->,thick] (X12) -- (X13) node [midway,above] {};
\draw [->,thick] (X13) -- (X14) node [midway,above] {};
\draw [->,thick] (X6) -- (X11) node [midway,above] {};
\draw [->,thick] (X7) -- (X12) node [midway,above] {};
\draw [double,-,thick] (X8) -- (X13) node [midway,above] {};
\draw [->,thick,dashed] (X7) -- (X11) node [midway,above] {};
\end{tikzpicture}
\end{center}
Since the top row is induced by an $n$-exact sequence, in particular $H(X^{n-1})\rightarrow H(X^n)$ is a weak cokernel of $H(X^{n-2})\rightarrow H(X^{n-1})$ in $\mathcal{M}$, and the morphism $H(X')\rightarrow F^n$ is a monomorphism, there exists a morphism $H(X^n)\dashrightarrow H(X')$ that makes the triangle containing $H(X^{n-1})$, $H(X^n)$ and $H(X')$ commutative. By \cite[Chapter 1, Proposition 5.6]{ARS} the left-hand square is a pushout square. The universal property of pushout implies that $H(X')\rightarrow F^n$ has a left inverse. Then the bottom row splits. This shows that $\Ext_{\mathcal{E}'}^1(\mathcal{M}, \mathcal{M})=0$.

Now assume that $k>1$. Using the argument of \cite[Remark 4.5]{EN}, we have the commutative diagram
\begin{diagram}\label{d1}
\begin{center}
\begin{tikzpicture}
\node (X1) at (-10,1) {$0$};
\node (X2) at (-8.5,1) {$H(X^0)$};
\node (X3) at (-7,1) {$\cdots$};
\node (X4) at (-5,1) {$H(X^{n-k})$};
\node (X5) at (-2,1) {$\cdots$};
\node (X6) at (0,1) {$H(X^{n-1})$};
\node (X7) at (2,1) {$H(X^n)$};
\node (X8) at (4,1) {$H(X^{n+1})$};
\node (X9) at (5.75,1) {$0$};
\node (X10) at (-6,-0.5) {$0$};
\node (X11) at (-4,-0.5) {$H(X')$};
\node (X12) at (-2,-0.5) {$\cdots$};
\node (X13) at (0,-0.5) {$F_{pb}^{n-1}$};
\node (X14) at (2,-0.5) {$H(X^n)$};
\node (X15) at (4,-0.5) {$H(X^{n+1})$};
\node (X16) at (5.75,-0.5) {$0$};
\draw [->,thick] (X1) -- (X2) node [midway,left] {};
\draw [->,thick] (X2) -- (X3) node [midway,left] {};
\draw [->,thick] (X3) -- (X4) node [midway,left] {};
\draw [->,thick] (X4) -- (X5) node [midway,left] {};
\draw [->,thick] (X5) -- (X6) node [midway,left] {};
\draw [->,thick] (X6) -- (X7) node [midway,above] {};
\draw [->,thick] (X7) -- (X8) node [midway,left] {};
\draw [->,thick] (X8) -- (X9) node [midway,above] {};
\draw [->,thick] (X10) -- (X11) node [midway,left] {};
\draw [->,thick] (X11) -- (X12) node [midway,above] {};
\draw [->,thick] (X12) -- (X13) node [midway,above] {};
\draw [->,thick] (X13) -- (X14) node [midway,above] {};
\draw [->,thick] (X14) -- (X15) node [midway,above] {};
\draw [->,thick] (X15) -- (X16) node [midway,above] {};
\draw [double,-,thick] (X7) -- (X14) node [midway,above] {};
\draw [double,-,thick] (X8) -- (X15) node [midway,above] {};
\end{tikzpicture}
\end{center}
\end{diagram}
with exact rows in $\mathcal{E}'$, where the bottom row is Yoneda equivalent to $\xi$.
Now for $j\in\{1,\cdots,n-1\}$ we set $C^j:=\Imm(H(X^j)\rightarrow H(X^{j+1}))$ in
$\mathcal{E}'$. Indeed, we split the first row of Diagram \eqref{d1} into short exact sequences as follows.
\begin{center}
\begin{tikzpicture}
\node (X1) at (-9,0) {$0$};
\node (X2) at (-7.5,0) {$H(X^0)$};
\node (X3) at (-5.5,0) {$H(X^1)$};
\node (X4) at (-3.5,0) {$H(X^2)$};
\node (X5) at (-1.5,0) {$\cdots$};
\node (X6) at (0.75,0) {$H(X^{n-1})$};
\node (X7) at (3,0) {$H(X^n)$};
\node (X8) at (5,0) {$H(X^{n+1})$};
\node (X9) at (6.75,0) {$0$};
\node (X10) at (-4.5,-1) {$C^1$};
\node (X11) at (-2.5,-1) {$C^2$};
\node (X12) at (-.5,-1) {$C^{n-2}$};
\node (X13) at (2,-1) {$C^{n-1}$};
\node (X14) at (-1.5,-1) {$\cdots$};
\draw [->,thick] (X1) -- (X2) node [midway,left] {};
\draw [->,thick] (X2) -- (X3) node [midway,left] {};
\draw [->,thick] (X3) -- (X4) node [midway,left] {};
\draw [->,thick] (X4) -- (X5) node [midway,left] {};
\draw [->,thick] (X5) -- (X6) node [midway,left] {};
\draw [->,thick] (X6) -- (X7) node [midway,above] {};
\draw [->,thick] (X7) -- (X8) node [midway,left] {};
\draw [->,thick] (X8) -- (X9) node [midway,above] {};
\draw [->>,thick] (X3) -- (X10) node [midway,left] {};
\draw [>->,thick] (X10) -- (X4) node [midway,above] {};
\draw [->>,thick] (X4) -- (X11) node [midway,above] {};
\draw [>->,thick] (X11) -- (X5) node [midway,above] {};
\draw [->>,thick] (X5) -- (X12) node [midway,above] {};
\draw [>->,thick] (X12) -- (X6) node [midway,above] {};
\draw [->>,thick] (X6) -- (X13) node [midway,above] {};
\draw [>->,thick] (X13) -- (X7) node [midway,above] {};
\end{tikzpicture}
\end{center}
First for $j\in\{1,\cdots,n-2\}$ by applying the functor $\Hom_{\mathcal{E}'}(-,H(X'))$ to the short exact sequence
\begin{equation}
0\rightarrow C^j\rightarrow H(X^{j+1})\rightarrow C^{j+1}\rightarrow 0, \notag
\end{equation}
we have the exact sequence\\
$
\Hom_{\mathcal{E}'}(H(X^{j+1}),H(X'))\rightarrow \Hom_{\mathcal{E}'}(C^j,H(X'))\rightarrow\Ext^1_{\mathcal{E}'}(C^{j+1},H(X'))\\\rightarrow \Ext^1_{\mathcal{E}'}(H(X^{j+1}),H(X'))=0. \notag
$
By the first part of the proof, $\Ext_{\mathcal{E}'}^1(\mathcal{M}, \mathcal{M})=0$ and hence $\Hom_{\mathcal{E}'}(H(X^{j+1}),H(X'))=0$.
Since $H(X^j)\rightarrow H(X^{j+1})$ is a weak cokernel of $H(X^{j-1})\rightarrow H(X^j)$ in $\mathcal{M}$ and $H(X^j)\rightarrow C^j$ is an epimorphism, $C^j\rightarrow H(X^{j+1})$ is a left $\mathcal{M}$-approximation, and so the map $\Hom_{\mathcal{E}'}(H(X^{j+1}),H(X'))\rightarrow \Hom_{\mathcal{E}'}(C^j,H(X'))$ is an epimorphism. Thus $\Ext_{\mathcal{E}'}^1(C^{j+1},H(X'))=0$ for $j\in\{1,\cdots,n-2\}$.

Now by induction assume that $k\in\{2,\cdots,n-1\}$ and $\mathcal{M}$ is $k$-rigid. We show that $\mathcal{M}$ is $(k+1)$-rigid.
Choose an arbitrary element $\xi\in \Ext_{\mathcal{E}'}^k(H(X^{n+1}),H(X'))$ and construct the Diagram \eqref{d1} for it. Since we have the following commutative diagram with exact row,
\begin{center}
\begin{tikzpicture}
\node (X1) at (-7.5,0) {$0$};
\node (X2) at (-6,0) {$H(X')$};
\node (X3) at (-4,0) {$\cdots$};
\node (X4) at (-2,0) {$F^{n-1}_{ab}$};
\node (X5) at (0,0) {$H(X^n)$};
\node (X6) at (2.5,0) {$H(X^{n+1})$};
\node (X7) at (4.5,0) {$0,$};
\node (X13) at (-1,-1) {$C^{n-1}$};
\draw [->,thick] (X1) -- (X2) node [midway,left] {};
\draw [->,thick] (X2) -- (X3) node [midway,left] {};
\draw [->,thick] (X3) -- (X4) node [midway,left] {};
\draw [->,thick] (X4) -- (X5) node [midway,left] {};
\draw [->,thick] (X5) -- (X6) node [midway,left] {};
\draw [->,thick] (X6) -- (X7) node [midway,above] {};
\draw [->>,thick] (X4) -- (X13) node [midway,above] {};
\draw [>->,thick] (X13) -- (X5) node [midway,above] {};
\end{tikzpicture}
\end{center}
it is enough to show that $\Ext_{\mathcal{E}'}^{k-1}(C^{n-1},H(X'))=0$. Because by hypothesis\\ $\Ext_{\mathcal{E}'}^{1,\cdots,k-1}(H(X^i),H(X'))=0$ for $i\in\{0,\cdots,n+1\}$, we have
$
\Ext_{\mathcal{E}'}^{k-1}(C^{n-1},H(X'))\cong \Ext_{\mathcal{E}'}^{k-2}(C^{n-2},H(X'))\cong \cdots \Ext_{\mathcal{E}'}^{2}(C^{n-k+2},H(X'))\cong \Ext_{\mathcal{E}'}^{1}(C^{n-k+1},H(X'))=0 \notag
$
and the result follows.
\end{proof}
\end{proposition}

\begin{theorem}\label{th1}
In the setting of Proposition \ref{pro10}, $\mathcal{M}=H(\mathcal{P})$ is an $n$-cluster tilting subcategory of $\mathcal{E}'=\frac{\mathcal{E}}{^{\bot}\mathcal{P}}$.
\begin{proof}
Denote by $\mathcal{Q}:\mathcal{E}\rightarrow \frac{\mathcal{E}}{^{\bot}\mathcal{P}}$ the quotient functor.
Let $Y=\mathcal{Q}(X)$ be an arbitrary object in $\mathcal{E}'$. By assumption there is a conflation
\begin{equation}\label{conf}
T\rightarrowtail X\twoheadrightarrow F
\end{equation}
such that $T\in^{\bot}\mathcal{P}$ and $F\in \cogen(\mathcal{P})$. By Lemma \ref{lem2}, there is a conflation $F\rightarrowtail \overline{F}\twoheadrightarrow T'$ where $T'\in^{\bot}\mathcal{P}$ and $\pd(\overline{F})\leq n-1$. Thus in the quotient category $\mathcal{E}'=\frac{\mathcal{E}}{^{\bot}\mathcal{P}}$ we have $\mathcal{Q}(X)\cong \mathcal{Q}(\overline{F})$ and there is an exact sequence
\begin{equation}
0\rightarrow P_n\rightarrow P_{n-1}\rightarrow \cdots\rightarrow P_1\rightarrow \overline{F}\rightarrow 0. \notag
\end{equation}
Applying $H$ to this exact sequence we obtain the following exact sequence in $\mathcal{E}'$
\begin{equation}
0\rightarrow H(P_n)\rightarrow H(P_{n-1})\rightarrow \cdots\rightarrow H(P_1)\rightarrow Y\rightarrow 0 \notag
\end{equation}
with $H(P_n),\cdots, H(P_1)\in \mathcal{M}$. Now because $\mathcal{P}$ and so $\mathcal{M}$ is closed under direct summands by \cite[Proposition 4.4]{Kv1} we have that $H(P_1)\twoheadrightarrow Y$ is a right $\mathcal{M}$-approximation and if $\Ext_{\mathcal{E}'}^i(H(\overline{F}),\mathcal{M})=0$ for $i\in\{1,\cdots,n-1\}$, then $Y\in \mathcal{M}$.

It remains to show that $\mathcal{M}$ is a cogenerating and covariantly finite subcategory. Consider the conflation \eqref{conf} for an arbitrary object $X$. Because $F\in \cogen(\mathcal{P})$ there is an inflation $i:F\rightarrowtail P$. In the category
$\mathcal{E}'$ we have $H(X)\cong H(F)$, and $H(i)$ is a inflation. This shows that $\mathcal{M}$ is a cogenerating subcategory.

Since by Proposition \ref{pro11}, $\mathcal{Q}:\mathcal{E}\rightarrow \frac{\mathcal{E}}{^{\bot}\mathcal{P}}$ preserve and reflect admissible morphisms, we have equivalence of exact categories
\begin{equation}
\mathcal{E}\simeq \modd_{\adm}(\mathcal{P})\simeq \modd_{\adm}(\mathcal{M}). \notag
\end{equation}
Because $\mathcal{M}$ is a generating contravariantly finite subcategory of $\mathcal{E}'$, by Proposition \ref{pro9}
\begin{equation}
\Upsilon_{\mathcal{M}}:\mathcal{E}'\longrightarrow \modd_{\adm}(\mathcal{M})\simeq \mathcal{E},\notag
\end{equation}
has a left adjoint $L_{\mathcal{M}}:\modd_{\adm}(\mathcal{M})\longrightarrow\mathcal{E}'$.
Now let $Y\in \mathcal{E}'$. By assumption there is a left $\mathcal{P}$-approximation $\Upsilon_{\mathcal{M}}(Y)\longrightarrow \Upsilon_{\mathcal{M}}(M)$ for some $M\in \mathcal{M}$. It can be easily seen that the corresponding map $Y\rightarrow M$ is a left $\mathcal{M}$-approximation for $Y$.
\end{proof}
\end{theorem}

Now we are ready to prove the main theorem.

\begin{theorem}\label{th2}
There is a bijection between the
following:
\begin{itemize}
\item[(1)]
Equivalence classes of $n$-cluster tilting subcategories of exact categories.
\item[(2)]
Equivalence classes of $n$-Auslander exact categories.
\end{itemize}
\begin{proof}
By Proposition \ref{pro12} and Proposition \ref{pro13}, for each $n$-cluster tilting subcategory $\mathcal{M}$ of an exact category $\mathcal{E}$ we have that $\modd_{\adm}(\mathcal{M})$ is an $n$-Auslander exact category. In the other word $\modd_{\adm}(-)$ defines a map from the first class to the second one.

Conversely let $\mathcal{E}$ be an $n$-Auslander exact category with $\mathcal{P}=\proj(\mathcal{E})$. By Theorem \ref{th1}, $\mathcal{E}'=\dfrac{\mathcal{E}}{^{\bot}\mathcal{P}}$ is an exact category and $\mathcal{E}'$ has an $n$-cluster tilting subcategory that is equivalent to $\mathcal{P}$. By Proposition \ref{pro14}, $\mathcal{E}\simeq\modd_{\adm}(\mathcal{P})\simeq \modd_{\adm}(\mathcal{M})$ and the result follows.
\end{proof}
\end{theorem}

\section{dominant dimension}
In this section we prove the more familiar version of higher Auslander correspondence for exact categories. Most of the proofs are similar to the \cite[Section 4.3]{HRK} and we omit them.

\begin{definition}$($\cite[Definition 4.11]{HRK}$)$
Let $\mathcal{E}$ be an exact category with enough projectives. The {\em dominant dimension} of $\mathcal{E}$,
denoted $\dom.dim(\mathcal{E})$, is the largest integer $n$ such that for any projective object $P\in \mathcal{E}$ there exists an
exact sequence $0\rightarrow P\rightarrow I_1\rightarrow \cdots\rightarrow I_n\rightarrow C\rightarrow0,$
with $I_k$ being projective and injective for $1\leq k\leq n$.
\end{definition}

\begin{theorem}\label{th3}
There is a bijection between the following:
\begin{itemize}
\item[(1)]
Equivalence classes of $n$-cluster tilting subcategories of exact categories with enough injectives.
\item[(2)]
Equivalence classes of exact categories $\mathcal{E}$ with enough projectives $\mathcal{P}=\proj(\mathcal{E})$ satisfying the following conditions:

$(a)$
$\gl.dim(\mathcal{E})\leq n+1\leq \dom.dim(\mathcal{E})$.

$(b)$
Any morphism $X\rightarrow E$ with $E\in ^{\bot}\mathcal{P}$ is admissible.

$(c)$
$\mathcal{P}$ is admissibly covariantly finite.
\end{itemize}
\begin{proof}
Let $\mathcal{M}$ be an $n$-cluster tilting subcategory of exact category $\mathcal{E}$ with enough injectives. By definition of $n$-cluster tilting subcategories, all injective objects belongs to $\mathcal{M}$. We show that $\modd_{\adm}(-)$ define a map from the first class to the seconde class. In order to show that $\modd_{\adm}(\mathcal{M})$ satisfies the conditions $(a)$, $(b)$ and $(c)$, it is enough to shows that $n+1\leq \dom.dim(\modd_{\adm}(\mathcal{M}))$. A projective object in $\modd_{\adm}(\mathcal{M})$ is isomorphic to $\Upsilon_{\mathcal{M}}(M)$ for some $M\in \mathcal{M}$. Choose an injective coresolution
\begin{equation}
0\rightarrow M\rightarrow I^0\rightarrow I^1\rightarrow \cdots\rightarrow I^n,\notag
\end{equation}
for $M$. Because $\mathcal{M}$ is $n$-rigid we have the following exact sequence in $\modd_{\adm}(\mathcal{M})$
\begin{equation}
0\rightarrow \Upsilon_{\mathcal{M}}(M)\rightarrow \Upsilon_{\mathcal{M}}(I^0)\rightarrow \Upsilon_{\mathcal{M}}(I^1)\rightarrow \cdots\rightarrow \Upsilon_{\mathcal{M}}(I^n)\rightarrow F\rightarrow 0,\notag
\end{equation}
where $F$ is the cokernel of $\Upsilon_{\mathcal{M}}(I^{n-1})\rightarrow \Upsilon_{\mathcal{M}}(I^n)$. Since $\gl.dim(\mathcal{E})\leq n+1$, $F$ is injective. Now by the proof of \cite[Lemma 4.10]{HRK}, $\Upsilon_{\mathcal{M}}(I^i)$ is projective and injective object for $i\in \{0,\cdots,n-1\}$, which proves the claim.

Conversely, let $\mathcal{E}$ be an exact category with enough projectives $\mathcal{P}=\proj(\mathcal{E})$ satisfying conditions $(a)$, $(b)$ and $(c)$. First we show that $\mathcal{E}$ is an $n$-Auslander exact category. We check the conditions of Definition \ref{def}.
\begin{itemize}
\item[(1)]
Because by assumption $\mathcal{P}$ is admissibly covariantly finite and $\dom.dim(\mathcal{E})\geq 1$, by \cite[Lemma 4.16]{HRK}, $(^{\bot}\mathcal{P},\cogen(\mathcal{P}))$ is a torsion pair.
\item[(2)]
Because $\dom.dim(\mathcal{E})\geq 1$, by \cite[Lemma 4.14]{HRK}, $^{\bot}\mathcal{P}$ is closed under admissible subobjects and hence by the condition $(b)$, $^{\bot}\mathcal{P}$ satisfies $\textbf{A2}$.
\item[(3)]
Because $\dom.dim(\mathcal{E})\geq n+1$, by the similar argument as in the proof of \cite[Lemma 4.13]{HRK}, we have $\Ext_{\mathcal{E}}^{1,\cdots,n}(^{\bot}\mathcal{P},\mathcal{P})=0$.
\item[(4)]
This is satisfied by assumption.
\item[(5)]
This is satisfied by assumption.
\end{itemize}
Now, the argument in the proof of $(2)\Rightarrow (1)$ of Theorem \ref{th2} and \cite[Lemma 4.12]{HRK} complete the proof.
\end{proof}
\end{theorem}

The following corollary is a higher Auslander correspondence for abelian categories \cite[Theorem 8.23]{Bel}

\begin{corollary}$($\cite[Theorem 8.23]{Bel}$)$
For any integer $n\geq0$, there is a bijection between the following:
\begin{itemize}
\item[(1)]
Equivalence classes of $n$-cluster tilting subcategories of abelian categories with enough injectives.
\item[(2)]
Equivalence classes of $n$-Auslander categories.
\end{itemize}
\begin{proof}
If $\mathcal{E}$ be an abelian category and $\mathcal{M}$ be an $n$-cluster tilting subcategory of it, then $\modd_{\adm}(\mathcal{M})=\modd(\mathcal{M})$ is an $n$-Auslander abelian category. Conversely if $\modd_{\adm}(\mathcal{M})$ be an abelian category, because $\eff(\mathcal{M})$ is a Serre subcategory of it, $\mathcal{E}\simeq\dfrac{\modd_{\adm}(\mathcal{M})}{\eff(\mathcal{M})}$ is an abelian category.
\end{proof}
\end{corollary}

\section*{acknowledgements}
The research of the first author was in part supported by a grant from IPM (No. 1400180047). Also, the research of the second author was in part supported by a grant from IPM (No. 1400170417).

\end{document}